\documentstyle[amscd, amstex]{amsart}

\def\ZZ         {{\Bbb Z}}
\def\RR         {{\Bbb R}}
\def\CC         {{\Bbb C}}

\def\PP         {{\Bbb P}}
\def\TT         {{\Bbb T}}

\def\rp         {{\cal RP}}
\def\gf        {{\cal GF}}

\newtheorem{thm}{Theorem}[section]
\newtheorem{lem}[thm]{Lemma}

\newtheorem{pr}[thm]{Proposition}

\theoremstyle{definition}
\newtheorem{rem}[thm]{Remark}

\newtheorem{defn}[thm]{Definition}

\newcommand{\pro}{{\rm Proj}}

\newcommand{\xds}{{{ X}_{\Delta^*}}}
\newcommand{\xd}{{{ X}_{\Delta}}}
\newcommand{\sd}{{ \Sigma_\Delta}}
\newcommand{\xsd}{{{X}_{ \Sigma_\Delta}}}

\newcommand{\co}{{\rm Conv}}

\newcommand{\im}{{\rm im}}
\newcommand{\spe}{{\rm Spec}}
\newcommand{\key}{\bibitem}

\newcommand{\xs}{{{ X}_{\Sigma}}}

\newcommand\hidot{{\raise1pt\hbox{$\scriptscriptstyle\bullet$}}}
\newcommand\lodot{{\raise.3pt\hbox{$\scriptscriptstyle\bullet$}}}

\newcommand{\Hom}{{\rm Hom}}

\newcommand{\ve}{{\scriptscriptstyle\vee}}
\begin{document}

\title[Degenerations and  mirror contractions]
{Degenerations and  mirror contractions \\ of Calabi-Yau complete
intersections  \\ via Batyrev-Borisov Mirror symmetry}
\author{Anvar R. Mavlyutov}
\address {Department of Mathematics, Oklahoma State  University, Stillwater, OK 74078, USA.}
 \email{mavlyutov@@math.okstate.edu}

%Fax: +1(405)-744-8275, mavlyutov@@math.okstate.edu \vskip1cm

% Math Subject Classification

\maketitle

\begin{abstract}
We show that the dual of the  Cayley  cone, associated to a
Minkowski sum decomposition of a reflexive polytope, contains a
reflexive polytope admitting a nef-partition. This nef-partition
corresponds to a Calabi-Yau complete intersection in a Gorenstein
Fano toric variety degenerating to an ample Calabi-Yau
hypersurface in another Fano toric variety. Using the
Batyrev-Borisov mirror symmetry construction, we found the mirror
contraction of a Calabi-Yau complete intersection to the mirror of
the ample Calabi-Yau hypersurface.
\end{abstract}

\tableofcontents
%\tableofcontents

\setcounter{section}{-1}

\section{Introduction.}

Refexive polytopes were introduced by V. Batyrev in \cite{b}.
These are lattice polytopes $\Delta$ in a real vector space
$\RR^d$ with lattice points in $\ZZ^d$ corresponding to monomials
of the anticanonical degree on a Gorenstein Fano toric variety.
Such polytopes are determined by the property that they have
vertices at lattice points and have the origin in their interior
 with the dual  polytope $\Delta^*=\{y\in  \RR^d\mid \langle \Delta,y\rangle\ge-1\}$    satisfying
the same property. This was the starting point for the Batyrev
construction of a large class of mirror pairs of Calabi-Yau
hypersurfaces in toric varieties in \cite{b}.

A Gorenstein Fano toric
 variety   associated to a reflexive polytope $\Delta$
 can be defined as $X_\Delta=\pro(\CC[\sigma\cap\ZZ^{d+1}])$, where $\sigma=\RR_{\ge0}\cdot(\Delta,1)\subset\RR^{d+1}$.
It contains an affine torus $\TT=(\CC^*)^d$ as a dense open subset
which acts naturally on the toric variety. A Calabi-Yau
hypersurface $Y_\Delta$ in a Gorenstein Fano toric
 variety  $X_\Delta$ can be viewed as the Zariski closure of a
 hypersurface
   $$\sum_{m\in\Delta\cap \ZZ^d} a_mt^m=0$$ in the  affine  torus
   $(\CC^*)^d\subset  X_\Delta$, where
   $m=(m_1,\dots,m_d)\in\ZZ^d$, $a_m\in\CC$
   and $t^m=t_1^{m_1}\cdots t_d^{m_d}$ for the coordinates $t_1,\dots,t_d$ on the torus.

More generally, a nef Calabi-Yau complete intersection in a
Gorenstein Fano toric variety $X_\Delta$ corresponds to a
Minkowski sum decomposition  of the reflexive polytope
 $\Delta=\Delta_0+\Delta_1+\cdots+\Delta_k$ by lattice polytopes.
The  Calabi-Yau complete intersection
$Y_{\Delta_0,\dots,\Delta_k}$ is the closure of the affine
complete intersection
$$\sum_{m\in\Delta_i\cap \ZZ^d} a_{i,m} t^m=0, \quad i=0,\dots,k$$
in  $(\CC^*)^d\subset  X_\Delta$ with generic coefficients $
a_{i,m}\in\CC$. A complete intersection in a toric variety is
called {\it nondegenerate} if every intersection  with a
$\TT$-orbit is either transversal or empty. A   generic nef
Calabi-Yau complete intersection is nondegenerate by Lemma~4.3 in
\cite{m0} and Proposition~6.8 in \cite{d}.

The mirror construction of Batyrev is a pair  of families of
nondegenerate Calabi-Yau hypersurfaces obtained as maximal
projective  crepant partial resolutions of $Y_\Delta$ and
$Y_{\Delta^*}$. Generalizing the polar duality of reflexive
polytopes, L.~Borisov in \cite{bo} introduced the notion of {\it
nef-partition}, which is a Minkowski sum decomposition
 of the reflexive polytope
 $\Delta=\Delta_0+\Delta_1+\cdots+\Delta_k$ by lattice polytopes
 such that the origin $0\in\Delta_i$ for all $i$.
A nef-partition has a dual nef-partition defined as the Minkowski
sum decomposition of the reflexive polytope
$\nabla=\nabla_0+\dots+\nabla_k$ in the dual vector space with
 $\nabla_j$   determined by  $\langle
\Delta_i,\nabla_j\rangle\ge-\delta_{ij}$ for all $0\le i,j\le k$,
where $\delta_{ij}$ is the Kronecker symbol. One of the basic
properties of the nef-partitions is that
$\Delta^*=\co(\nabla_0,\dots,\nabla_k)$ and
$\nabla^*=\co(\Delta_0,\dots,\Delta_k)$. The Batyrev-Borisov
mirror symmetry construction is a pair of families of
nondegenerate nef Calabi-Yau complete intersections obtained as
maximal projective  crepant partial resolutions of
$Y_{\Delta_0,\dots,\Delta_k}$ and $Y_{\nabla_0,\dots,\nabla_k}$.

A topological mirror symmetry test for compact $n$-dimensional
Calabi-Yau manifolds $V$ and $V^*$ is a symmetry of their Hodge
numbers: $h^{p,q}(V)=h^{n-p,q}(V^*), 0\le p,q\le n.$ For singular
varieties Hodge numbers must be replaced by the stringy Hodge
numbers $h^{p,q}_{\rm st}$ introduced by V. Batyrev in \cite{b2}.
The usual Hodge numbers coincide with the stringy Hodge numbers
for nonsingular Calabi-Yau varieties. Moreover, all crepant
partial resolutions $\widehat V$ of singular Calabi-Yau varieties
$V$ have the same stringy Hodge numbers: $h^{p,q}_{\rm
st}(\widehat V)=h^{p,q}_{\rm st}(  V )$. In \cite{bb3}, Batyrev
and Borisov show that the pair of Calabi-Yau complete
intersections $V= Y_{\Delta_0,\dots,\Delta_k}$  and
$V^*=Y_{\nabla_0,\dots,\nabla_k}$ pass the mirror symmetry test.
One of the main ingredients of their proof was the use of the
Cayley trick which associates to a Calabi-Yau complete
intersection $Y_{\Delta_0,\dots,\Delta_k}$ a generalized
Calabi-Yau hypersurface in a higher dimensional Fano toric variety
$\pro(\CC[\bar\sigma\cap \ZZ^{d+k+1} ])$, where
 $\bar\sigma= \{ ( \sum_{i=0}^k
t_i\Delta_i,t_0,\dots,t_k)\mid t_i\in\RR_{\ge0}\}\subset
\RR^{d+k+1}$, called the {\it Cayley cone} associated to the
polytopes $\Delta_0,\dots,\Delta_k$.

There are six different reflexive polytopes of dimension $d+k$
associated to a Minkowski sum decomposition of a $d$-dimensional
reflexive polytope $\Delta=\Delta_0+\Delta_1+\cdots+\Delta_k$ into
$k+1$ lattice polytopes. One of them, contained in the Cayley cone
$\bar\sigma$ at an integral distance $k+1$ from the origin, is
isomorphic to
\begin{equation}\label{e:ref1} (k+1){\rm Conv}(
  \Delta_0, \Delta_1+e_1,\dots,\Delta_k +e_k)- \sum_{i=1}^k
  e_i\end{equation}
in $\RR^{d+k}\simeq\RR^d\oplus\RR^k$ where $\{e_1,\dots,e_k\}$ is
the standard basis of  $\RR^k$.  The dual $\bar\sigma^\ve$ of the
Cayley cone also contains a reflexive polytope isomorphic to
\begin{equation}\label{e:ref2} (k+1)\co\biggl(\biggl\{u-\sum_{i=1}^k \min\langle
\Delta_i,u\rangle
 e_i^*\mid u\in \Delta^*\biggr\}\cup\{e_1^*,\dots,e_k^*\}\biggl)- \sum_{i=1}^k e_i^*\end{equation} in
 $\RR^{d+k}\simeq\RR^d\oplus\RR^{k}$, where $\{e_1^*,\dots,e_k^*\}$ is
the standard basis of  $\RR^k$ dual to  $\{e_1,\dots,e_k\}$. The
reflexive polytopes (\ref{e:ref1}) and  (\ref{e:ref2}) are not
dual to each other, and their dual polytopes give another two
reflexive polytopes. It turns out that while the polytope
(\ref{e:ref1}) may not admit a nef-partition (if
$\Delta=\Delta_0+\Delta_1+\cdots+\Delta_k$ is not a
nef-partition), the reflexive polytope (\ref{e:ref2}) always
admits one: $\hat\nabla_0+\dots+\hat\nabla_k$, where
$\hat\nabla_0=\co(\{u-\sum_{i=1}^k \min\langle \Delta_i,u\rangle
 e_i^*\mid u\in \Delta^*\}\cup\{e_1^*,\dots,e_k^*\})$,
  $\hat\nabla_i=\hat\nabla_0-e_i^*$, for $i=1,\dots,k$.
The dual of this nef-partition is
$\tilde\Delta_0+\cdots+\tilde\Delta_k$, where
$\tilde\Delta_i=\co(\Delta_i+e_i,0)$, for $i=0\dots,k$, and
$e_0:=-\sum_{i=1}^ke_i$. This  reflexive polytope together with
its dual $\co(\hat\nabla_0,\dots,\hat\nabla_k)$ are the other two
reflexive polytopes associated to a Minkowski sum decomposition of
the reflexive polytope $\Delta$.

The dual of $\hat\nabla_0+\dots+\hat\nabla_k$ is the reflexive
polytope $\co(\tilde\Delta_0,\dots,\tilde\Delta_k)$. The
Gorenstein Fano toric variety
$X_{\hat\nabla_0+\dots+\hat\nabla_k}$, whose fan consists of the
 cones over the proper faces of
 $\co(\tilde\Delta_0,\dots,\tilde\Delta_k)$,
 is the ambient space of   deformations of the Gorenstein Fano toric
 variety
$X_{\Delta^*}$ in \cite{m3}. In this paper, we show that the
embedding $X_{\Delta^*}\hookrightarrow
X_{\hat\nabla_0+\cdots+\hat\nabla_k}$ realizes the ample
Calabi-Yau hypersurface $Y_{\Delta^*}$ as a  complete intersection
in $X_{\hat\nabla_0+\dots+\hat\nabla_k}$, which deforms to a
nondegenerate Calabi-Yau complete intersection
$Y_{\hat\nabla_0,\dots,\hat\nabla_k}$ corresponding to the
nef-partition $\hat\nabla_0+\dots+\hat\nabla_k$. The degeneration
$Y_{\hat\nabla_0,\dots,\hat\nabla_k}\rightsquigarrow Y_{\Delta^*}$
can be lifted to the degeneration of a maximal  projective crepant
partial resolution $Y'_{\hat\nabla_0,\dots,\hat\nabla_k}$ of
$Y_{\hat\nabla_0,\dots,\hat\nabla_k}$ to a partial resolution
$Y'_{\Delta^*}$ of $Y_{\Delta^*}$. Taking maximal projective
crepant partial resolution   $Y''_{\Delta^*}$ of $Y'_{\Delta^*}$
we obtain a geometric transition (a contraction followed by
smoothing) from a minimal Calabi-Yau hypersurface to a minimal
Calabi-Yau complete intersection:
$Y'_{\hat\nabla_0,\dots,\hat\nabla_k}\rightsquigarrow
Y'_{\Delta^*}\leftarrow Y''_{\Delta^*}$. According to a conjecture
of D.~Morrison in \cite{mor}, every geometric transition between
Calabi-Yau manifolds should correspond to a mirror geometric
transition between the mirror partners of the original Calabi-Yau
manifolds with the roles of degeneration and contraction reversed.
In Section~\ref{s:mirror}, we explicitly construct a natural
contraction of a minimal Calabi-Yau complete intersection
$Y'_{\tilde\Delta_0,\dots,\tilde\Delta_k}$ to a degenerate
Calabi-Yau hypersurface $ Y'_{\Delta}$ in a maximal projective
crepant partial resolution of $X_{\Delta}$. The smoothing of $
Y'_{\Delta}$ to a nondegenerate Calabi-Yau hypersurface $
Y''_{\Delta}$ gives a geometric transition
$Y''_\Delta\rightsquigarrow Y'_{\Delta}\leftarrow
Y'_{\tilde\Delta_0,\dots,\tilde\Delta_k}$, which should be the
mirror of the above one. We use the method of   Batyrev in
\cite{b1}, \cite{bs} to support the mirror correspondence of
geometric transitions by showing that the degeneration of  the
hypergeometric series arising from the main period of   Calabi-Yau
varieties coincides with the hypergeometric series of the maximal
projective partial crepant resolution of the degenerate
Calabi-Yau. These hypergeometric series determine the mirror map
between the K\"ahler and complex moduli spaces (see
\cite[Sec.~6.3.4]{ck}).

The construction of deformations of ample Calabi-Yau hypersurfaces
and their partial resolutions are consistent with our conjecture
in \cite{m2} that all deformations of Calabi-Yau complete
intersections (of dimension $\ge3$) in toric varieties are
Calabi-Yau complete intersections in higher dimensional toric
varieties. An application of deformations of Gorenstein Fano toric
varieties to deformations of nef Calabi-Yau complete intersections
and a generalization of the above geometric transitions will
appear in \cite{m6}. These constructions together with the
previously known geometric transitions between Calabi-Yau
hypersurfaces in \cite{bkk,mor} give a strong evidence that the
web of Calabi-Yau complete intersections in toric varieties can be
connected by explicit geometric transitions.

Here is an organization of our paper. In Section~\ref{s:refl}, we
study properties of reflexive Gorenstein cones and explicitly
describe the Cayley cone and its dual together with the reflexive
polytopes contained in theses cones. Then we briefly overview some
basic notation and facts of toric geometry. Section~\ref{s:caydef}
explains the relation of the Cayley trick and deformations of Fano
toric varieties constructed in \cite{m3}, and
Section~\ref{s:defcy} constructs deformations of  Calabi-Yau
hypersurfaces. Finally, in Section~\ref{s:mirror}, we construct
two geometric transitions described above, and then
Section~\ref{s:degen} discusses degenerations of the main periods
of Calabi-Yau complete intersections and Mirror Symmetry.

{\it Acknowledgment.} We would like to thank Victor Batyrev
 for  pointing out the reference \cite{bn}.

\section{Combinatorics of reflexive polytopes and Gorenstein cones.} \label{s:refl}

In this section, we explicitly describe the reflexive polytopes
arising from the construction of the Cayley cone associated to a
Minkowski sum decomposition of a reflexive polytope  in \cite{bb}.
We show that the dual of the Cayley cone contains a reflexive
polytope which   admits a nef-partition introduced in \cite{bo}.

 Let $N$ be a lattice and $M$ be its dual lattice with a paring $\langle *,*\rangle:M\times N\rightarrow\ZZ$.

\begin{defn}\cite{b}  A lattice polytope $\Delta$ in
$M_\RR=M\otimes\RR$  (i.e., its vertices are at the lattice
points) is called a {\it reflexive polytope} if it contains $0$ in
its interior and   the  dual polytope
$$\Delta^*=\{n\in N_\RR\mid \langle m,n\rangle\ge-1\, \forall\, m\in \Delta\}$$
in the dual vector space $N_\RR=N\otimes \RR$ is also a lattice
polytope. The pair $\Delta$ and $\Delta^*$ is called a pair of
dual reflexive polytopes and it satisfies $\Delta=(\Delta^*)^*$.
\end{defn}

Reflexive polytopes are related to the  notion of reflexive
Gorenstein cones from \cite{bb}. Let $\bar M$ and $\bar N$ be
lattices   which are dual to each other.  Let $\sigma\subset \bar
M_\RR$ be a  polyhedral cone with a vertex at $0$. The dual cone
of $\sigma$ is defined as $$\sigma^\ve=\{n\in \bar N_\RR\mid
\langle m, n\rangle\ge0 \,\forall\,  m\in\sigma \}.$$

\begin{defn}\cite{bb} A maximal dimensional  polyhedral  cone
 $\sigma$ is
called {\em Gorenstein}, if  it is generated by finitely many
lattice points contained in the affine hyperplane  $\{x\in \bar
M\mid
 \langle x,n_\sigma\rangle=1\}$ for a unique $n_{\sigma}\in \bar N$. A Gorenstein cone $\sigma$ is called   {\it reflexive}
 if both $\sigma$ and $\sigma^\ve$ are Gorenstein
cones, in which case they both have maximal dimension and uniquely
determined $n_\sigma\in  \bar N$ and $m_{\sigma^\ve}\in \bar M$,
which take value 1 at the primitive lattice generators of the
respective cones. The positive integer $r=\langle
m_{\sigma^\ve},n_\sigma\rangle$ is called the index of the
reflexive Gorenstein cones $\sigma$ and $\sigma^\ve$.
\end{defn}

Denote   $\sigma_{(i)}:=\{x\in\sigma\mid\langle
x,n_\sigma\rangle=i\}$,  the slice of the cone at an integral
distance $i$ from the origin. Since  $r=\langle
m_{\sigma^\ve},n_\sigma\rangle$,  the lattice points
$m_{\sigma^\ve}$ and $n_\sigma $ lie in the interiors of
$\sigma_{(r)}=r\cdot\sigma_{(1)}$ and
$\sigma^\ve_{(r)}=r\cdot\sigma^\ve_{(1)}$, respectively.

\begin{pr}\cite[Pr.~2.11]{bb} Let $\sigma$ be a Gorenstein cone. Then
$\sigma$ is a reflexive Gorenstein cone of index $r$ if and only
if the polytope $\sigma_{(r)}-m_{\sigma^\ve}$ is a reflexive
polytope with respect to the lattice $\bar M\cap
n_\sigma^\perp=\{x\in \bar M\mid \langle x,n_\sigma\rangle=0\}$.
\end{pr}

  As noted in Remark 1.13  in \cite{bn}, the
  reflexive
polytopes $\sigma_{(r)}-m_{\sigma^\ve}$ and
$\sigma^\ve_{(r)}-n_\sigma$   are combinatorially dual to each
other, but not dual as lattice polytopes. In \cite[Proposition
1.15]{bn}, the  dual reflexive polytope
$(\sigma_{(r)}-m_{\sigma^\ve})^*$ was obtained as
$\sigma^\ve_{(1)}$ with respect to the refined affine lattice
$(\bar N+\frac{1}{r}\ZZ n_\sigma)\cap  \{y\in \bar N_\RR\mid
 \langle m_{\sigma^\ve}, y\rangle=1\}$.
We will give an alternative description for the dual reflexive
polytope.

\begin{pr}\label{p:gorefl} Let $\sigma$     be a  reflexive Gorenstein cone  of index
$r$.
 Then the dual polytope
 $(\sigma_{(r)}-m_{\sigma^\ve})^*=\pi\bigl(\sigma^\ve_{(1)}\bigr)$
with respect to the lattice  $\bar N /\ZZ n_\sigma\simeq \Hom(\bar
M\cap n_\sigma^\perp,\ZZ)$,
 where $\pi:\bar N_\RR\rightarrow \bar N_\RR/\RR n_\sigma$   is  the  quotient homomorphism.
\end{pr}

\begin{pf}   Note that the vertices $v$  of
$\sigma^\ve_{(1)}$ are in one-to-one correspondence with the
facets  $F_v:=\sigma_{(r)}\cap v^\perp$  of the polytope
 $\sigma_{(r)}$ by the duality of the cones $\sigma$ and $\sigma^\ve$. Then $\bigl\langle
F_v-m_{\sigma^\ve}, v \bigr\rangle=-1$ with respect to the pairing
of $\bar M$ and $\bar N$. Consequently,   $\bigl\langle
F_v-m_{\sigma^\ve}, \pi( v) \bigr\rangle=-1$ with respect to the
pairing of $\bar M\cap n_\sigma^\perp$ and $\bar N /\ZZ n_\sigma$
Hence, all vertices of the dual polytope
$(\sigma_{(r)}-m_{\sigma^\ve})^*$ in $\bar N_\RR/\RR n_\sigma$ are
of the form $\pi(v)$ for  a vertex   $v$ of $\sigma^\ve_{(1)}$.
\end{pf}

 A special class of reflexive Gorenstein cones arises from a Calabi-Yau complete intersection in
  a Gorenstein Fano toric variety by a Cayley trick (see \cite{bb}).
  Let
$\Delta$ be a reflexive polytope in   $M_\RR$ and
$\Delta=\Delta_0+\Delta_1+\cdots+\Delta_k$ be a Minkowski sum
decomposition   by lattice polytopes. By
\cite[Proposition~3.6]{bb}, the cone
$$\bar\sigma=\Biggl\{\biggl( \sum_{i=0}^k
t_i\Delta_i,t_0,\dots,t_k\biggr) \mid
t_i\in\RR_{\ge0}\Biggr\}\subset M_\RR\oplus\RR^{k+1}
$$ is reflexive Gorenstein of index $k+1$. This cone is called the
{\it Cayley cone} associated to the polytopes
$\Delta_0,\dots,\Delta_k$. It can also be written as
$$\bar\sigma=\RR_{\ge0}\cdot{\rm Conv}(
  \Delta_0+r_0, \Delta_1+r_1,\dots,\Delta_k +r_k),$$
where $\{r_0,\dots,r_k\}\subset \ZZ^{k+1}\subset M\oplus\ZZ^{k+1}$
is the standard
 basis of the second summand. A third way to write the Cayley
 cone is given in the following lemma.

\begin{lem}\label{l:2eq} There is equality of cones  $$\bar\sigma =\bigl\{\bigl(
t \cdot{\rm Conv}(
  \Delta_0, \Delta_1+e_1,\dots,\Delta_k +e_k),t\bigr)\mid t\in\RR_{\ge0}\bigr\},$$
 induced by the isomorphism
$$M\oplus\ZZ^{k+1}\simeq M \oplus\ZZ^{k }\oplus\ZZ,\quad
 (m,\alpha_0,\dots,\alpha_k)  \mapsto (m,\alpha_1,\dots,\alpha_k,
\alpha_0+\cdots+\alpha_k),$$ where $\{e_1,\dots,e_k\}$ is the
standard basis for the second summand $\ZZ^{k }$.
\end{lem}

 The dual of the Cayley cone $\bar\sigma^\ve$ can also be
 explicitly found.

\begin{pr} \label{p:cdual}  Let
$\bar\sigma\subset M_\RR\oplus\RR^{k+1}$ be the Cayley cone
associated to $\Delta_0,\dots,\Delta_k$. Then
$$\bar\sigma^\ve=\RR_{\ge0}\cdot\co\biggl(\biggl\{u-\sum_{i=0}^k \min\langle
 \Delta_i,u\rangle
 r^*_i\mid u\in \Delta^*\biggr\}\cup\{r^*_0, \dots,r^*_k\}\biggr),$$
 where $\{r^*_0,\dots,r_k^*\}$ is the   basis of $\ZZ^{k+1}\subset N\oplus\ZZ^{k+1}$ dual to
 $\{r_0,\dots,r_k\}$.
\end{pr}

\begin{pf} We have $u+\sum_{i=0}^k\alpha_i r_i^*\in\bar\sigma^\ve$
with $u\in N_\RR$ and $\alpha_i\in\RR$ if and only if $\langle
x_j+r_j, u+\sum_{i=0}^k\alpha_i r_i^*\rangle\ge0$ for all
$x_j\in\Delta_j$, $j=0,\dots,k$. But the last inequality is
equivalent to $\alpha_j\ge-\langle x_j, u\rangle$ for all
$x_j\in\Delta_j$. Hence, $\alpha_j\ge-\min\langle
\Delta_j,u\rangle$. Since $0$ is in the interior of $\Delta$,
$\min\langle \Delta,u\rangle<0$ for $u\ne0$, whence
$u+\sum_{i=0}^k\alpha_i r_i^*= -\min\langle
 \Delta ,u \rangle(u'-\sum_{i=0}^k \min\langle
 \Delta_i,u'\rangle
 r^*_i)+\sum_{i=0}^k\beta_i r_i^*$, where $u'=-\frac{u}{\min\langle
 \Delta,u\rangle}\in\Delta^*$ and $\beta_i=\alpha_i+\min\langle
 \Delta_i,u\rangle\ge0$.
\end{pf}

The following alternative  view of the dual of the Cayley cone may
also be useful.

\begin{lem}\label{l:1eq}
There is  equality of cones
 $$\bar\sigma^\ve =\biggl\{\Bigl(t\cdot\co\Bigl(\Bigl\{u-\sum_{i=1}^k \min\langle
 \Delta_i,u\rangle e_i^*\mid u\in \Delta^*\Bigr\}\cup\{e_1^*,\dots,e_k^*\}\Bigr),t\Bigr)\mid t\in\RR_{\ge0}\biggr\}$$ induced by the
isomorphism
$$ N \oplus\ZZ^{k+1}\simeq N \oplus\ZZ^{k }\oplus\ZZ,\quad
 (n,\alpha_0,\dots,\alpha_k)  \mapsto (n,\alpha_1,\dots,\alpha_k,
\alpha_0+\cdots+\alpha_k),$$
 where $\{e_1^*,\dots,e_k^*\}$ is the
standard basis for the second summand $\ZZ^{k }$.
\end{lem}

The above descriptions of the Cayley cone $\bar\sigma$, associated
to a Minkowski sum decomposition of the reflexive polytope
$\Delta=\Delta_0+\Delta_1+\cdots+\Delta_k$,  and  of the dual cone
$\bar\sigma^\ve$ directly show that both cones are Gorenstein
reflexive of index $k+1$ with the unique lattice points
$n_{\bar\sigma}=r_0^*+r_1^*+\cdots+r_k^*$ and
$m_{\bar\sigma^\ve}=r_0+r_1+\cdots+r_k$. Next, we want to
explicitly describe the reflexive polytopes
$\bar\sigma_{(k+1)}-m_{{\bar\sigma}^\ve}$ and
$\bar\sigma^\ve_{(k+1)}-n_{{\bar\sigma}}$, and their dual
polytopes, arising from Proposition~\ref{p:gorefl}.

Consider the lattice $\tilde{M}:=M\oplus\ZZ^k$ and denote by
$\{e_1,\dots,e_k\}$ the standard basis for the second summand
 $\ZZ^k$.  Then $\tilde{N}:=N\oplus\ZZ^k$ is the dual to $\tilde{M}$
 lattice and set
 $\{e_1^*,\dots,e_k^*\}$ be the dual to $\{e_1,\dots,e_k\}$ basis
 in $\ZZ^k$.
The following statements follow  trivially from
Propositions~\ref{p:gorefl} and \ref{p:cdual}.

\begin{lem}\label{l:1st}  Let
$\bar\sigma\subset M_\RR\oplus\RR^{k+1}$ be the Cayley cone
associated  to lattice polytopes $\Delta_0,\dots,\Delta_k$ in
$M_\RR$ such that $\Delta=\Delta_0+\Delta_1+\cdots+\Delta_k$ is
reflexive. Then
 $$\bar\sigma_{(k+1)}-m_{{\bar\sigma}^\ve}
\simeq(k+1){\rm Conv}(
  \Delta_0, \Delta_1+e_1,\dots,\Delta_k +e_k)- \sum_{i=1}^k e_i$$  induced by the
isomorphism   $\bar M\cap  n_{\bar\sigma}^\perp\simeq \tilde M$,
 $
 (m,\alpha_0,\dots,\alpha_k)  \mapsto
(m,\alpha_1,\dots,\alpha_k),$  where $\bar M:=M\oplus\ZZ^{k+1}$
and $n_{\bar\sigma}=\sum_{i=0}^k r_i^*$.  Also,
$$\pi\bigl(\bar\sigma_{(1)}\bigr)={\rm Conv}\biggl(
  \Delta_0-\sum_{i=1}^k e_i, \Delta_1+e_1,\dots,\Delta_k +e_k\biggr),$$
where
 $\pi:\bar M\rightarrow
\bar M /\ZZ m_{{\bar\sigma}^\ve}\simeq\tilde M$,
 $(m,\alpha_0,\dots,\alpha_k)  \mapsto
(m,\alpha_1-\alpha_0,\dots,\alpha_k-\alpha_0).$
\end{lem}

\begin{lem}  \label{l:2nd} Let
$\bar\sigma\subset M_\RR\oplus\RR^{k+1}$ be the Cayley cone
associated  to lattice polytopes $\Delta_0,\dots,\Delta_k$ in
$M_\RR$ such that $\Delta=\Delta_0+\Delta_1+\cdots+\Delta_k$ is
reflexive.   Then
 $$\bar\sigma^\ve_{(k+1)}-n_{{\bar\sigma}}
\simeq(k+1)\co\biggl(\Bigl\{u-\sum_{i=1}^k \min\langle
 \Delta_i,u\rangle e_i^*\mid u\in \Delta^*\Bigr\}\cup\{e_1^*,\dots,e_k^*\}\biggr)- \sum_{i=1}^k e_i^*$$   under the
isomorphism  $\bar N\cap m_{{\bar\sigma}^\ve}^\perp\simeq \tilde
N$,
 $(n,\alpha_0,\dots,\alpha_k)  \mapsto
(n,\alpha_1,\dots,\alpha_k),$  where $\bar N:=N \oplus\ZZ^{k+1}$,
$m_{{\bar\sigma}^\ve}=\sum_{i=0}^k r_i$. Also,
$$\pi\bigl(\bar\sigma^\ve_{(1)}\bigr)=\co\biggl(\biggl\{u-\sum_{i=0}^k
\min\langle
 \Delta_i,u\rangle
 e_i^*\mid u\in \Delta^*\biggr\}\cup\{e_0^*, \dots,e_k^*\}\biggr),$$ where
 $\pi:\bar N\rightarrow
 \bar N/\ZZ n_{\bar\sigma}\simeq\tilde N$,
 $(n,\alpha_0,\dots,\alpha_k)  \mapsto
(n,\alpha_1-\alpha_0,\dots,\alpha_k-\alpha_0)$, and
$e_0^*:=-\sum_{i=1}^k e_i^*$.
\end{lem}

It is straightforward to check  that the above  natural
isomorphisms $\bar M\cap n_{\bar\sigma}^\perp\simeq \tilde M$,
$\bar N/\ZZ n_{\bar\sigma}\simeq\tilde N $ and $\bar N\cap
m_{{\bar\sigma}^\ve}^\perp\simeq \tilde N$, $\bar M /\ZZ
m_{{\bar\sigma}^\ve}\simeq\tilde M$ respect pairings. Hence, by
combining the above lemmas with Proposition~\ref{p:gorefl} we get
two pairs of reflexive polytopes.

\begin{pr}\label{p:2iso} Let $\Delta$ be a reflexive polytope in $M_\RR$ and
$\Delta=\Delta_0+\Delta_1+\cdots+\Delta_k$ be a Minkowski sum
decomposition by lattice polytopes in $M_\RR$. Denote $\widehat
\Delta_0= {\rm Conv}(
  \Delta_0, \Delta_1+e_1,\dots,\Delta_k +e_k)$,
$\widehat\Delta_i=\widehat\Delta_0-e_i$ in $\tilde M_\RR$,  for
$i=1,\dots,k$. Then $\widehat\Delta_0+\cdots+\widehat\Delta_k
=(k+1)\widehat\Delta_0+e_0$ is a reflexive polytope in
$\tilde{M}_\RR=M_\RR\oplus\RR^k$ with the dual reflexive polytope
$$
(\widehat\Delta_0+\cdots+\widehat\Delta_k)^*=
\co\biggl(\biggl\{u-\sum_{i=0}^k \min\langle
 \Delta_i,u\rangle
 e_i^*\mid u\in \Delta^*\biggr\}\cup\{e_0^*, \dots,e_k^*\}\biggr)$$
  in
$\tilde{N}_\RR=N_\RR\oplus\RR^k$, where
  $e_0=-\sum_{i=1}^k e_i$, $e_0^*=-\sum_{i=1}^k e_i^*$.
\end{pr}

\begin{pr}\label{p:2ndrefl}  Let $\Delta$ be a reflexive polytope in $M_\RR$ and
$\Delta=\Delta_0+\Delta_1+\cdots+\Delta_k$ be a Minkowski sum
decomposition by lattice polytopes in $M_\RR$. Denote
$\hat\nabla_0=\co(\{u-\sum_{i=1}^k \min\langle  \Delta_i,u\rangle
 e_i^*\mid u\in \Delta^*\}\cup\{e_1^*,\dots,e_k^*\}),
  \hat\nabla_i=\hat\nabla_0-e_i^*$ in $\tilde{N}_\RR$, for $i=1,\dots,k$. Then
 $\hat\nabla_0+\dots+\hat\nabla_k=(k+1)\hat\nabla_0+e_0^*$ is a reflexive polytope in
$\tilde{N}_\RR=N_\RR\oplus\RR^k$ with the dual reflexive polytope
$$ (\hat\nabla_0+\dots+\hat\nabla_k)^*={\rm Conv}(
  \Delta_0+e_0, \Delta_1+e_1,\dots,\Delta_k +e_k)
 $$
  in
$\tilde{M}_\RR=M_\RR\oplus\RR^k$, where  $e_0=-\sum_{i=1}^k e_i$,
 $e_0^*=-\sum_{i=1}^k e_i^*$.
\end{pr}

It turns out that one of the reflexive polytopes in the above
propositions always admits  a nef-partition introduced in
\cite{bo}. A {\it nef-partition} of a reflexive polytope $\Delta$
is a Minkowski sum decomposition $\Delta=\Delta_0+\dots+\Delta_k$
by lattice polytopes such that the origin $0\in\Delta_i$ for all
$i$.
  If one defines the polytopes $$\nabla_j=\{y\in\RR^d\mid\langle
 x,y\rangle\ge-\delta_{ij}\, \forall \,x\in\Delta_i,\, i=0,\dots,k\}$$ for $j=0,\dots,k$, then
   $\nabla=\nabla_0+\dots+\nabla_k$ is a reflexive polytope and  $\nabla_j$ are lattice polytopes with $0\in \nabla_j$
for all $j$.
 The nef-partitions $\Delta=\Delta_0+\dots+\Delta_k$ and
 $\nabla=\nabla_0+\dots+\nabla_k$ are called dual to each other.

For convenience, we will introduce the following notation. For any
two subsets $P$ and $Q$ of a real vector space we denote by
$P\uplus Q:=\co(P\cup Q)$, the convex hull of the union of $P$ and
$Q$. The operation $\uplus$ is clearly associative and
commutative.

By Propositions 3.1 and 3.2 in \cite{bo}, we have the following
dualities.

\begin{pr}\cite{bo}\label{p:bo} Let $\Delta$ be a reflexive polytope in $M_\RR$
and let $\Delta=\Delta_0+ \cdots+\Delta_k$ be a nef-partition and
$\nabla=\nabla_0+ \cdots+\nabla_k$  be the dual nef-partition in
$N_\RR$. Then
 $( \Delta_0\uplus\cdots\uplus \Delta_k)^*=\nabla_0+\cdots+\nabla_k$ and $(\nabla_0\uplus\cdots\uplus \nabla_k)^*=\Delta_0+\cdots+ \Delta_k.
$
\end{pr}

Now, suppose $\Delta=\Delta_0+\Delta_1+\cdots+\Delta_k$ is a
Minkowski sum decomposition of a reflexive polytope by lattice
polytopes in $M_\RR$. By construction, we see that
$\hat\nabla_0+\dots+\hat\nabla_k$ is a nef-partition even if
$\Delta=\Delta_0+ \cdots+\Delta_k$ is not a nef-partition, since
$0\in \hat\nabla_i$, for $i=0,\dots,k$, and the polytopes  $
\hat\nabla_i$ are convex hulls of lattice points.  It is not
difficult to determine that the dual nef-partition is
$\tilde\Delta_0+\cdots+\tilde\Delta_k$, where
$\tilde\Delta_i:=\co(\Delta_i+e_i, 0 )$. In particular, we have

\begin{pr}\label{p:1iso} Let $\Delta$ be a reflexive polytope in $M_\RR$ and
$\Delta=\Delta_0+ \cdots+\Delta_k$ be a Minkowski sum
decomposition   by lattice polytopes in $M_\RR$. Then
$\hat\nabla_0+\dots+\hat\nabla_k$ is a reflexive polytope in
$\tilde{N}_\RR=N_\RR\oplus\RR^k$. Moreover,
$\hat\nabla_0+\dots+\hat\nabla_k$ is a nef-partition dual to
$\tilde\Delta_0+\cdots+\tilde\Delta_k$  and the following
identities hold: $$
(\hat\nabla_0\uplus\dots\uplus\hat\nabla_k)^*=\tilde\Delta_0+\cdots+\tilde\Delta_k,\quad
   (\tilde\Delta_0\uplus\cdots\uplus\tilde\Delta_k)^*=\hat\nabla_0+\dots+\hat\nabla_k.
  $$
\end{pr}

 One can easily find the lattice points in
$\tilde\Delta_0\uplus\cdots\uplus\tilde\Delta_k$ and
$\hat\nabla_0\uplus\dots\uplus\hat\nabla_k$. Denote by $L(P)$ and
$l(P)$
  the set  and the number of lattice points in a polytope $P$ in a real
vector space.

\begin{pr}\label{p:lat}  Let $\Delta$ be a reflexive polytope in $M_\RR$ and
$\Delta=\Delta_0+ \cdots+\Delta_k$ be a Minkowski sum
decomposition   by lattice polytopes in $M_\RR$. Then\\
$L(\tilde\Delta_0\uplus\cdots\uplus\tilde\Delta_k)=\{0\}
 \cup\bigcup_{i=0}^kL(\Delta_i+e_i) $,\quad $l(\tilde\Delta_0\uplus\cdots\uplus\tilde\Delta_k)=1+\sum_{i=0}^k
l(\Delta_i) $,  $L(\hat\nabla_0)=  \{n-\sum_{i=1}^k \min\langle
\Delta_i,n\rangle
 e_i^*\mid n\in \Delta^*\cap N \}\cup\{e_1^*,\dots,e_k^*\}$,\quad
 $l(\hat\nabla_0)=l(\Delta^*)+k$,
 $l(\hat\nabla_0\uplus\dots\uplus\hat\nabla_k)=
(k+1)l(\Delta^*)+k^2 $.   \end{pr}

The   reflexive polytope
$\tilde\Delta_0\uplus\cdots\uplus\tilde\Delta_k$ arises from the
construction of deformations of Gorenstein Fano toric varieties
$X_{\Delta^*}$ associated to the fan generated by the faces of
$\Delta=\Delta_0+\cdots+\Delta_k$ in \cite{m3}. The deformations
are realized by complete intersections in a higher dimensional
Fano toric variety whose fan is generated by the faces of
$\tilde\Delta_0\uplus\cdots\uplus\tilde\Delta_k$. On the other
hand, the reflexive polytope
$\widehat\Delta_0+\cdots+\widehat\Delta_k$ arises from the Cayley
trick: it corresponds to the dual of the canonical line bundle (or
the anticanonical degree) on the projective space bundle
associated to a Calabi-Yau complete intersection in the Fano toric
variety $X_\Delta$ whose fan is generated by the faces of
$\Delta^*$. We will review these constructions in detail in
Section~\ref{s:caydef}.

From Lemma~\ref{l:1st}, we can see that the reflexive polytopes $
\tilde\Delta_0\uplus\cdots\uplus\tilde\Delta_k $ and $
\widehat\Delta_0+\cdots+\widehat\Delta_k $ must be the same up to
a linear transformation and a change of the lattice. The same
should hold for the dual reflexive polytopes. These
transformations can be explicitly described as follows.

\begin{lem}\label{l:isoref1} The homomorphism of lattices $$\varphi: M \oplus\ZZ^{k }\rightarrow  M \oplus\ZZ^{k
},\quad  m+\sum_{i=1}^k\alpha_i e_i   \mapsto
 (k+1)m+\sum_{i=1}^k\alpha_i ((k+1)e_i+e_0),$$
 maps $\tilde\Delta_0\uplus\cdots\uplus\tilde\Delta_k$ onto
$\widehat\Delta_0+\cdots+\widehat\Delta_k$.
\end{lem}

\begin{lem}\label{l:isoref2} The homomorphism of lattices $$\varphi^*: N \oplus\ZZ^{k }\rightarrow  N \oplus\ZZ^{k
},\quad  n+\sum_{i=1}^k\alpha_i e_i^*   \mapsto
 (k+1)n+\sum_{i=1}^k\alpha_i ((k+1)e_i^*+e_0^*),$$
 maps $(\widehat\Delta_0+\cdots+\widehat\Delta_k)^*$  onto
$(\tilde\Delta_0\uplus\cdots\uplus\tilde\Delta_k)^*$.
\end{lem}

Finishing this section, we will look at what happens if
$\Delta=\Delta_0 +\cdots+\Delta_k$ is a nef-partition. In this
case the dual nef-partition $\nabla=\nabla_0+\dots+\nabla_k$
 satisfies $\langle
\Delta_i,\nabla_j\rangle\ge-\delta_{ij}$ for all $0\le i,j\le k$,
and the dual to $\Delta$ reflexive polytope is $\Delta^* ={\rm
Conv}(
  \nabla_0,\dots,\nabla_k  )$ with $\nabla_i\cap\nabla_j=\{0\}$ for all $i,j$. Now if
   $0\ne u\in\nabla_i\cap N$, then  $-1=\min\langle
 \Delta,u\rangle=\min\langle \Delta_0,u\rangle+\cdots+\min\langle \Delta_k,u\rangle$.
Since $\langle \Delta_j,u\rangle\ge0$ for $j\ne i$ and $\langle
\Delta_i,u\rangle\ge-1$, we conclude that $\min\langle
\Delta_j,u\rangle=0$ for $j\ne i$ and $\min\langle
\Delta_i,u\rangle=-1$. Hence, by Proposition~\ref{p:cdual}, we get
$$\bar\sigma^\ve=\RR_{\ge0}\cdot{\rm Conv}(
  \nabla_0+r_0^*,  \dots,\nabla_k +r_k^*).$$
Similarly, $\hat\nabla_0={\rm Conv}(
  \nabla_0, \nabla_1+e_1^*,\dots,\nabla_k +e_k^*)=\widehat\nabla_0$. Applying Proposition~\ref{p:bo}, we get the following dualities
  for eight reflexive polytopes of dimension $d+k$ corresponding to a dual pair of
  nef-partitions of dimension $d$.

\begin{pr}\label{p:ref}  Let $\Delta$ be a reflexive polytope in $M_\RR$
and let $\Delta=\Delta_0+ \cdots+\Delta_k$ be a nef-partition and
$\nabla=\nabla_0+ \cdots+\nabla_k$  be the dual nef-partition in
$N_\RR$. Then
$$(\tilde\Delta_0\uplus\cdots\uplus\tilde\Delta_k)^*=\widehat\nabla_0+\cdots+\widehat\nabla_k\text{ and }
(\tilde\nabla_0\uplus\cdots\uplus\tilde\nabla_k)^*=\widehat\Delta_0+\cdots+\widehat\Delta_k
$$ are nef-partitions  respectively dual to
$$ (\widehat\nabla_0\uplus\cdots\uplus\widehat\nabla_k)^* =\tilde\Delta_0+\cdots+\tilde\Delta_k\text{ and }
 (\widehat\Delta_0\uplus\cdots\uplus\widehat\Delta_k)^*=
 \tilde\nabla_0+\cdots+\tilde\nabla_k.
$$
\end{pr}

\section{Some basics from toric geometry.}\label{s:basic}

This section will review some basic facts from  \cite{c},
 \cite{c2}, \cite{f}    on toric geometry.
See  \cite{d}, \cite{o} for additional references.

 Let $\xs$ be a $d$-dimensional toric variety associated
with a finite rational polyhedral fan $\Sigma$ in $N_\RR$. Denote
by $\Sigma(1)$ the finite set of the 1-dimensional cones $\rho$ in
$\Sigma$, which correspond to the torus invariant divisors
$D_\rho$ in $\xs$.    From the work of David Cox (see \cite{c}),
every toric variety can be described as a categorical quotient of
a Zariski open subset of an affine space by a subgroup of a torus.
For simplicity, assume that the   1-dimensional cones  $\Sigma(1)$
span $N_\RR$. Consider the polynomial ring
$S(\Sigma):=\CC[x_\rho\mid \rho\in\Sigma(1)]$, called the {\it
homogeneous coordinate ring} of the toric variety $\xs$, and the
corresponding affine space $\CC^{\Sigma(1)}=\spe(\CC[x_\rho\mid
\rho\in\Sigma(1)])$. Let $B =\langle
\prod_{\rho\not\subseteq\sigma} x_\rho\mid \sigma\in\Sigma\rangle$
be the ideal in $S(\Sigma)]$. This ideal determines a Zariski
closed set ${\bf V}(B )$ in $\CC^{\Sigma(1)}$, which
 is
invariant under the diagonal group action of the subgroup
$$G=\biggl\{(\mu_\rho) \in(\CC^*)^{\Sigma(1)}\mid\prod_{\rho\in\Sigma(1)} \mu_\rho^{\langle
u,v_\rho \rangle}=1\, \forall\, u\in M\biggr\}$$ of the torus
$(\CC^*)^{\Sigma(1)}$ on the affine space $\CC^{\Sigma(1)}$, where
$v_\rho$ denotes the primitive lattice generator of the
1-dimensional cone $\rho$. Then by Theorem 2.1 in \cite{c}, the
toric variety $\xs$ is the categorical quotient
$(\CC^{\Sigma(1)}\setminus{\bf V}(B))/G$. This presentation is
important because it allows us to work with closed subvarieties of
the toric variety. In particular, a torus invariant divisor
$D_\rho$ is given by the equation $x_\rho=0$.

The ring $S(\Sigma)$ is graded by the the Chow group
$$A_{d-1}(\xs)\simeq \Hom(G,\CC^*),$$
 and  $\deg(\prod_{\rho\in\Sigma(1) }
x_\rho^{b_\rho})= [\sum_{\rho\in\Sigma(1)} b_\rho D_\rho]\in
A_{d-1}(\xs)$. For a torus invariant Weil divisor
$D=\sum_{\rho\in\Sigma(1)} b_\rho D_\rho$, there is a one-to-one
correspondence between the monomials of $\CC[x_\rho:
\rho\in\Sigma(1)]$ in the degree $[\sum_{\rho\in\Sigma(1)} b_\rho
D_\rho]\in A_{d-1}(\xs)$ and the lattice points inside the
polytope
$$ \Delta_D=\{m\in M_\RR\mid\langle m,v_\rho\rangle\ge-b_\rho\,\forall\,
\rho\in\Sigma(1)\}
$$ by associating  to $m\in\Delta_D$ the monomial  $ \prod_{\rho\in\Sigma(1) }
x_\rho^{b_\rho+\langle m,v_\rho\rangle}=x^m\prod_{\rho\in\Sigma(1)
} x_\rho^{b_\rho}$ where $x^m$ will denote $
\prod_{\rho\in\Sigma(1) } x_\rho^{\langle m,v_\rho\rangle}$. If we
denote the  homogeneous degree of $S(\Sigma)$ corresponding to
$\beta=[D]\in A_{d-1}(\xs)$ by $S(\Sigma)_\beta$, then by
Proposition 1.1 in \cite{c}, we also have a natural isomorphism
$$H^0(\xs,O_{\xs}(D))\simeq S(\Sigma)_\beta.$$
In particular, every hypersurface in $\xs$ of degree
$\beta=\sum_{\rho\in\Sigma(1)} b_\rho D_\rho$ corresponds to a
polynomial $$\sum_{m\in\Delta_D\cap M} a_m \prod_{\rho\in\Sigma(1)
} x_\rho^{b_\rho+\langle m,v_\rho\rangle}$$ with the coefficients
$a_m\in \CC$.

Every lattice polytope $\Delta$ in $M_\RR$ determines the Weil
divisor $$D_\Delta=\sum_{\rho\in\Sigma(1)} -\min\langle
\Delta,v_\rho\rangle D_\rho$$ on $\xs$.  By Theorem 1.6 in
\cite{m01} we know that if $D$ is a Cartier divisor on a compact
toric variety $\xs$, then $O_\xs(D)$ is generated by global
sections iff $D$ is numerically effective (nef).   In this case,
by \cite[p.~68]{f}, we
 get  $D=D_{\Delta_D}$. Also, if
  for a lattice
polytope $\Delta$ the divisor $D_\Delta$ is nef, then
 $\Delta_{D_\Delta}=\Delta$.  Additionally, this correspondence preserves sums:   if
$D_{\Delta_1}$ and $D_{\Delta_2}$ are nef then
$D_{\Delta_1+\Delta_2}=D_{\Delta_1}+D_{\Delta_2}$. Moreover, the
 following holds.

\begin{lem}\label{l:nef} Let $\xs$ be a compact toric variety   associated to a fan $\Sigma$ in
$N_\RR$. Suppose ${\Delta_1}$ and $ {\Delta_2}$ are lattice
polytopes in $M_\RR$ then $D_{\Delta_1+\Delta_2}$ is a nef divisor
on $\xs$ iff $D_{\Delta_1}$ and $D_{\Delta_2}$ are nef on $\xs$.
\end{lem}

\begin{pf} If  $D_{\Delta_1+\Delta_2}$ is a nef divisor  on $\xs$ then $\Sigma$ is a refinement of the normal fan of
$\Delta_1+\Delta_2$. But since $\Delta_i$ is a Minkowski summand
of $ {\Delta_1+\Delta_2}$, the normal fan of $\Delta_1+\Delta_2$
is a refinement of the normal fans of $\Delta_i$, for $i=1,2$.
Hence, $\Sigma$ is a refinement of the normal fans of $\Delta_1$
and $\Delta_2$. This implies that $D_{\Delta_1}$ and
$D_{\Delta_2}$ are nef on $\xs$, if $D_{\Delta_1+\Delta_2}$ is a
nef divisor on $\xs$. The other direction follows from the fact
that the sum of nef divisors is nef.
\end{pf}

   From  Mori's theory we know that nef divisors
correspond to contractions and for toric varieties this
correspondence can be formulated as in Theorem 1.2 in \cite{m1}.

\begin{thm}\label{t:fun}  Let $[D]\in A_{d-1}(\xs)$ be a nef divisor class on
a compact toric variety $\xs$ of dimension $d$. Then, there exists
a unique compact toric variety $X_{\Sigma_D}$ with a surjective
toric morphism $\pi:\xs@>>>X_{\Sigma_D}$ such that $\pi^*[Y]=[D]$
for some ample divisor $Y$ on $X_{\Sigma_D}$. Moreover, $\dim
X_{\Sigma_D}=\dim \Delta_D$, and the fan
$\Sigma_D=\Sigma_{\Delta_D}$, the normal fan of polytope
$\Delta_D$, for a torus invariant $D$.
\end{thm}

Finishing this section, we will recall an alternative way to
describe projective toric varieties using the language of
Gorenstein cones. Suppose that $\Delta$ is a lattice polytope in
$M_\RR$ such that its support function
$\psi_\Delta=-\min\langle\Delta,\underline{\,\,\,}\,\rangle$ is
strictly convex with respect to the fan $\Sigma$. In this case,
the divisor $D_\Delta$ is ample and $\Sigma=\sd$ is the normal fan
of $\Delta$. Consider the Gorenstein  cone $$K=\{( t\Delta,t)\mid
t\in\RR_{\ge0}\}\subset M_\RR\oplus\RR.$$
 The projective toric variety
 $\xd:=\xsd$
can be represented as $ {\rm Proj}({\Bbb C}[K\cap(M\oplus\ZZ)]).$
 Moreover, if  $\beta\in A_{d-1}(\xd)$ is
the class of the  ample divisor
$D_\Delta=\sum_{\rho\in\Sigma_\Delta(1)} b_\rho D_\rho$, then
there is a natural isomorphism of graded rings
\begin{equation}\label{e:isom}
{\Bbb C}[K\cap(M\oplus\ZZ)] \simeq\bigoplus_{i=0}^\infty
S(\sd)_{i\beta},
\end{equation}
sending $\chi^{(m,i)}\in\CC[K\cap(M\oplus\ZZ)]_i$ to
$\prod_{\rho\in\sd(1)} x_\rho^{i b_\rho+\langle
m,v_\rho\rangle}=x^m\prod_{\rho\in\sd(1)} x_\rho^{i b_\rho}$. This
correspondence allows to translate an equation  of a hypersurface
given by a polynomial  in homogeneous coordinates
$$\sum_{m\in\Delta\cap M}a_m \prod_{\rho\in\sd(1)} x_\rho^{i
b_\rho+\langle m,v_\rho\rangle}= \sum_{m\in\Delta\cap M}a_m x^m
\prod_{\rho\in\sd(1)} x_\rho^{i b_\rho}
$$ into the homogeneous  element $\sum_{m\in\Delta\cap M}a_m \chi^{(m,i)}$
 of the
Gorenstein ring.

\section{Cayley trick and deformations of Fano toric varieites.}\label{s:caydef}

To describe the Cayley trick used in mirror symmetry by \cite{bb}
we start with    a Gorenstein Fano toric variety $X_\Delta:
=\xsd$, whose (normal) fan $\Sigma_\Delta$
 of the reflexive polytope $\Delta$ consists of  the cones   generated by the proper faces of the dual reflexive polytope
$\Delta^*$ in $N_\RR$. Consider a Minkowski sum decomposition
$\Delta=\Delta_0+ \cdots+\Delta_k$ by lattice polytopes. The
anticanonical divisor $D_\Delta=\sum_{\rho\in\Sigma_\Delta(1)}
D_\rho$ on the Fano toric variety $\xd$ is ample, and, in
particular, nef. Applying Lemma~\ref{l:nef}, we get the nef
divisors $D_{\Delta_0}, \dots,D_{\Delta_k}$ on $\xd$. Given a
collection of line bundles on a variety, the Cayley trick
associates to it the projective space bundle.  In our case we get
the $\PP^{k}$-bundle $\PP({\cal E}_{\Delta_0,
\dots,\Delta_k})\rightarrow\xd,$ where
$${\cal E}_{\Delta_0, \dots,\Delta_k}={\cal
O}_\xd(D_{\Delta_0})\oplus  \cdots\oplus {\cal
O}_\xd(D_{\Delta_k}).$$ By \cite[p.~58]{o}, we know that this
bundle is a toric variety with its fan in $N_\RR\oplus \RR^k$.

\begin{pr} The torus invariant anticanonical divisor on  $\PP({\cal
E}_{\Delta_0,\dots,\Delta_k})$ is big and nef and equals
$D_{\widehat\Delta_0+\cdots+\widehat\Delta_k}$.
\end{pr}

\begin{pf} We only need to check that for the torus invariant anticanonical divisor $Y$
of the toric variety
 $\PP({\cal
E}_{\Delta_0,\dots,\Delta_k})$ there is equality of polytopes:
$\Delta_Y=\widehat\Delta_0+\cdots+\widehat\Delta_k$ in
$M_\RR\oplus \RR^k$. But this follows immediately from the fan
description in \cite[p.~58]{o} and Proposition~\ref{p:2iso}.
\end{pf}

By Theorem~\ref{t:fun}, for the nef divisor
$D_{\widehat\Delta_0+\cdots+\widehat\Delta_k}$  we get the
contraction $$\PP({\cal E}_{\Delta_0,\dots,\Delta_k})\rightarrow
X_{\widehat\Delta_0+\cdots+\widehat\Delta_k},$$ which relates the
projective  bundle to the Fano toric  variety. In the case, when
the anticanonical divisor on  $\PP({\cal
E}_{\Delta_0,\dots,\Delta_k})$  is ample (i.e., the vector bundle
${\cal E}_{\Delta_0, \dots,\Delta_k}$ is ample) we get the
equality $\PP({\cal E}_{\Delta_0,\dots,\Delta_k})=
X_{\widehat\Delta_0+\cdots+\widehat\Delta_k}$.  The projective
bundle $ \PP({\cal E}_{\Delta_0,\dots,\Delta_k})\rightarrow\xd$
can also be viewed as a   contraction corresponding to the
polytope $\Delta$ in $M_\RR\subset M_\RR\oplus \RR^k$ and its nef
divisor on  $ \PP({\cal E}_{\Delta_0,\dots,\Delta_k})$.

The Fano toric varieties $\xd$ and
$X_{\widehat\Delta_0+\cdots+\widehat\Delta_k}$ can also be
described in the language of Gorenstein cones from \cite{bb}. Let
$\sigma=\{( t\Delta,t)\mid t\in\RR_{\ge0}\}\subset M_\RR\oplus\RR$
and
$$\bar\sigma=\biggl\{\biggl( \sum_{i=0}^k
t_i\Delta_i,t_0,\dots,t_k\biggr)\mid
t_i\in\RR_{\ge0}\biggr\}\subset M_\RR\oplus\RR^{k+1}.$$ Then, by
the correspondence at the end of Section~\ref{s:basic} and
Lemma~\ref{l:2eq}, we have $\xd=\pro(\CC[\sigma\cap
(M\oplus\ZZ)])$ and
$X_{\widehat\Delta_0+\cdots+\widehat\Delta_k}=\pro(\CC[\bar\sigma\cap
\bar M ])$, where $\bar M=M\oplus\ZZ^{k+1}$. Inclusion of cones
$\sigma\subset\bar\sigma$ via
$$M_\RR\oplus\RR\hookrightarrow M_\RR\oplus\RR^{k+1},\quad
(m,r)\mapsto(m,r,\dots,r),$$ induces an injective homomorphism
$\CC[\sigma\cap (M\oplus\ZZ)]\hookrightarrow\CC[\bar\sigma\cap
\bar M]$ and
 a surjective morphism
 $\spe(\CC[\bar\sigma\cap  \bar M])\rightarrow
\spe(\CC[\sigma\cap (M\oplus\ZZ)])$ of affine toric varieties. It
also induces a rational map $\pro(\CC[\bar\sigma\cap \bar M])
\dashrightarrow \pro(\CC[\sigma\cap (M\oplus\ZZ)])$ of projective
toric varieties.  This map coincides with the morphism $ \PP({\cal
E}_{\Delta_0,\dots,\Delta_k})\rightarrow\xd$, if ${\cal
E}_{\Delta_0, \dots,\Delta_k}$ is an ample vector bundle.

There is more story to the Cayley trick in associating a semiample
hypersurface in the projective bundle to the nef Calabi-Yau
complete intersection on $\xd$ given by global sections of ${\cal
O}_\xd(D_{\Delta_0}),\dots,{\cal O}_\xd(D_{\Delta_k})$, but we
will not need this here.

 Now, let us show how the Cayley trick is related to deformations of Fano toric varieties.   Consider the Fano toric variety
 $\xds$,  whose
 fan $\Sigma_{\Delta^*}$ in $M_\RR$ consists of  the cones   generated by the proper faces of the  reflexive polytope
$\Delta=(\Delta^*)^*$. Take the same Minkowski sum decomposition
$\Delta=\Delta_0+ \dots+\Delta_k$ as above. We have a natural
inclusion of spaces $M_\RR\subset M_\RR\oplus \RR^k$ which induces
the inclusion of  polytopes $ \Delta\subset
(k+1)(\tilde\Delta_0\uplus\cdots\uplus\tilde\Delta_k)$ and the map
of fans over the proper faces of these polytopes.

\begin{thm}\cite{m3} Associated to the map of fan
$\Sigma_{\Delta^*}$ to
$\Sigma_{(\tilde\Delta_0\uplus\cdots\uplus\tilde\Delta_k)^*}$, the
  toric morphism $X_{\Delta^*}\rightarrow
 X_{ {(\tilde\Delta_0\uplus\cdots\uplus\tilde\Delta_k)^*}}$ is an embedding, whose image  is a complete
intersection given by the equations $$\prod_{
 v_\rho\in\Delta_i+e_i } x_\rho -\prod_{   v_\rho\in \Delta_0+e_0} x_\rho
=0,$$ for $i=1,\dots,k$, where $x_\rho$ are the homogeneous
coordinates of the toric variety $X_{
{(\tilde\Delta_0\uplus\cdots\uplus\tilde\Delta_k)^*}}$
corresponding to the vertices $v_\rho$ of the polytope $
\tilde\Delta_0\uplus\cdots\uplus\tilde\Delta_k $.
\end{thm}

Let $l(\Delta^*)$ denotes the number of lattice points in the
reflexive polytope $\Delta^*$.   By \cite{m3}, we have $
(kl(\Delta^*)-k) $-parameter embedded deformation family of $\xds$
in $X_{ {(\tilde\Delta_0\uplus\cdots\uplus\tilde\Delta_k)^*}}$
given by the equations:
$$\Biggl(x^{ e_i^*}-1+
 \sum_{n\in(\Delta^*\cap N)\setminus\{0\} }\lambda_{i,n}
x^{  n-\sum_{j=1}^k \min\langle \Delta_j,n\rangle e_j^*
}\Biggr)\prod_{
 v_\rho\in\Delta_0+e_0 } x_\rho=0 $$ for
$i=1,\dots,k$.

The embedding $X_{\Delta^*}\hookrightarrow
 X_{ {(\tilde\Delta_0\uplus\cdots\uplus\tilde\Delta_k)^*}}$ can also be
 described in the language of Gorenstein cones. Let $\sigma$ and
 $\bar\sigma$ be the same cones as above. Associated to the
 inclusion of cones
$\sigma\subset\bar\sigma$, there is a projection
$\bar\sigma^\ve\rightarrow\sigma^\ve$ induced by $$\bar
N:=N\oplus\ZZ^{k+1}\rightarrow N\oplus\ZZ,\quad
(n,\alpha_0,\dots,\alpha_k)\mapsto(n,\alpha_0+\cdots+\alpha_k)$$
and the corresponding  ring homomorphism
$$ \CC[\bar\sigma^\ve\cap \bar N]\longrightarrow\CC[\sigma^\ve\cap
(N\oplus\ZZ)],$$  which is surjective by Lemma 2.2 in \cite{m3}.
Hence, we get the embedding $$\spe(\CC[\sigma^\ve\cap
(N\oplus\ZZ)])\hookrightarrow\spe(\CC[\bar\sigma^\ve\cap \bar N])
$$ of affine toric varieties. By (\ref{e:isom}) and Lemma~\ref{l:1eq}, if $\beta=\deg(\prod_{
 v_\rho\in\Delta_0+e_0 } x_\rho)=[D_{\hat\nabla_0}]$, then
$$\CC[\bar\sigma^\ve\cap \bar N]\simeq\bigoplus_{i=0}^\infty
S(\Sigma_{\hat\nabla_0})_{i\beta},\quad \chi^{u+\sum_{j=0}^k
\alpha_j r^*_j  }\mapsto x^{u+\sum_{j=1}^k \alpha_j e^*_j }\prod_{
 v_\rho\in\Delta_0+e_0 } x_\rho^{\alpha_0+\cdots+\alpha_k}.$$  Since
$\Sigma_{\hat\nabla_0}=\Sigma_{(\tilde\Delta_0\uplus\cdots\uplus\tilde\Delta_k)^*}$
by Proposition~\ref{p:1iso},  we also get the embedding of
projective toric varieties
$$X_{\Delta^*}=\pro(\CC[\sigma^\ve\cap (N\oplus\ZZ)])
\hookrightarrow X_{
{(\tilde\Delta_0\uplus\cdots\uplus\tilde\Delta_k)^*}}=\pro(\CC[\bar\sigma^\ve\cap
\bar N]),$$ where the image is a complete intersection given by
$\chi^{r_i^*}-\chi^{r_0^*}$,   for $i=1,\dots,k$.

 Then deformations of the Fano toric variety $X_{\Delta^*}$
are $\pro(\CC[\bar\sigma^\ve\cap \bar N ]/I)$, where the ideal
$I\subset\CC[\bar\sigma^\ve\cap \bar N] $ is generated by
$$ \chi^{r_i^*}-\chi^{r_0^*}+\sum_{n \in(\Delta^*\cap N)\setminus\{0\} }\lambda_{i,n}\chi^{
n-\sum_{j=0}^k \min\langle \Delta_j,n\rangle r_j^* },$$
 for $i=1,\dots,k$, where  $\{r_0^*,\dots,r_k^*\}$ is the
 basis of  $\ZZ^{k+1}\subset N\oplus\ZZ^{k+1}$.

The ambient toric variety $X_{
{(\tilde\Delta_0\uplus\cdots\uplus\tilde\Delta_k)^*}}=\pro(\CC[\bar\sigma^\ve\cap
\bar N])$ of the deformation of $X_{\Delta^*}$ is related to the
Fano toric variety
$X_{\widehat\Delta_0+\cdots+\widehat\Delta_k}=\pro(\CC[\bar\sigma\cap
\bar M ])$ from the Cayley trick by the duality of the Gorenstein
cones. Note that the reflexive polytopes associated to these toric
varieties are not dual to each other, but a precise relation
between them is described in Lemmas~\ref{l:isoref1} and
\ref{l:isoref2}.

We will conclude this section by  considering   the case when
$\Delta=\Delta_0+\dots+\Delta_k$ is a nef-partition in $M_\RR$ and
$\nabla=\nabla_0+ \cdots+\nabla_k$ is the dual nef-partition in
$N_\RR$. In this case, by Proposition~\ref{p:ref}, we get
 $X_{(\tilde\Delta_0\uplus\cdots\uplus\tilde\Delta_k)^*}=X_{\widehat\nabla_0+\cdots+\widehat\nabla_k}$
and
$X_{(\tilde\nabla_0\uplus\cdots\uplus\tilde\nabla_k)^*}=X_{\widehat\Delta_0+\cdots+\widehat\Delta_k}$.
The fan of the projective space bundle $\PP({\cal
 E}_{\nabla_0,\dots,\nabla_k})$ is a refinement of the normal fan of $\widehat\nabla_0+\cdots+\widehat\nabla_k$,
which is obtained by a subdivision of the faces of the reflexive
polytope $\tilde\Delta_0\uplus\cdots\uplus\tilde\Delta_k
 =(\widehat\nabla_0+\cdots+\widehat\nabla_k)^*$. Intersection of
 this fan with the linear subspace $M_\RR\subset\tilde  M_\RR\oplus
 \RR^k$ gives a subdivision $\Sigma'_{\Delta^*}$ of the normal
 fan $\Sigma_{\Delta^*}$ of the polytope $\Delta^*$.

 By the
the embeddings of   toric varieties from \cite[Section~7]{m3} we
have a commutative diagram:
\begin{equation*}
\begin{array}{ccccc}
X_{\Sigma'_{\Delta^*}}&\hookrightarrow & \PP({\cal
 E}_{\nabla_0,\dots,\nabla_k})&\rightarrow &  X_\nabla\\
 \downarrow&&\downarrow& \\
X_{\Delta^*}&\hookrightarrow &
 X_{ {(\tilde\Delta_0\uplus\cdots\uplus\tilde\Delta_k)^*}}&=&
 X_{\widehat\nabla_0+\cdots+\widehat\nabla_k}
\end{array}
\end{equation*}
Similarly, if  $\Sigma'_{\nabla^*}$ is obtained by intersecting
the fan of $\PP({\cal
 E}_{\Delta_0,\dots,\Delta_k})$ with the subspace $N_\RR\subset\tilde  N_\RR\oplus
 \RR^k$, then
\begin{equation*}
\begin{array}{ccccc}
X_{ \Sigma'_{\nabla^*}}&\hookrightarrow & \PP({\cal
 E}_{\Delta_0,\dots,\Delta_k})&\rightarrow &  X_\Delta\\
 \downarrow&&\downarrow& \\
X_{\nabla^*}&\hookrightarrow &
 X_{ {(\tilde\nabla_0\uplus\cdots\uplus\tilde\nabla_k)^*}}&=&
 X_{\widehat\Delta_0+\cdots+\widehat\Delta_k}.
\end{array}
\end{equation*}

\section{Deformations of Calabi-Yau hypersurfaces in Fano toric
varieties.}\label{s:defcy}

In this section we   show that deformations of Fano toric
varieties induce deformations of Calabi-Yau hypersurfaces. The
embedding of the ambient Fano toric variety realizes a Calabi-Yau
hypersurface as a Calabi-Yau complete intersection in a higher
dimensional Fano toric variety. The deformations of the resulting
complete intersections are ``polynomial'', corresponding to
changing the  coefficients at the monomials. As before, we assume
for the rest that $\Delta$ is a reflexive polytope and
$\Delta=\Delta_0+\Delta_1+\cdots+\Delta_k$ is a Minkowski sum
decomposition by lattice polytopes in $M_\RR$.

\begin{thm}  Let
$Y_{\Delta^*}\subset X_{\Delta^*}=\pro(\CC[\sigma^\ve\cap
(N\oplus\ZZ)])$ be an ample Calabi-Yau hypersurface given by the
equation $$ \sum_{n\in\Delta^*\cap N}a_n \chi^{(n,1)}
 =0,$$ where $a_n\in\CC$.
Then the image of $ Y_{\Delta^*} $ by the embedding $
X_{\Delta^*}\hookrightarrow
 X_{ {(\tilde\Delta_0\uplus\cdots\uplus\tilde\Delta_k)^*}}=\pro(\CC[\bar\sigma^\ve\cap
\bar N])$   is a nef Calabi-Yau complete intersection given by the
equations $$
 a_0\chi^{r_0^*}+\sum_{n \in(\Delta^*\cap N)\setminus\{0\} }a_n \chi^{
n-\sum_{j=0}^k \min\langle \Delta_j,n\rangle r_j^* }  =0,\quad
 \chi^{r_i^*}-\chi^{r_0^*}=0, \, i=1,\dots,k.$$

\end{thm}

\begin{pf}  We need
to show that the kernel of the surjective $\ZZ$-graded ring
homomorphism $$ \CC[\bar\sigma^\ve\cap \bar
N]\longrightarrow\CC[\sigma^\ve\cap (N\oplus\ZZ)]/(f),$$ where $f=
\sum_{n\in\Delta^*\cap N}a_n \chi^{(n,1)}$, is generated by $$\bar
f=a_0\chi^{r_0^*}+\sum_{n \in(\Delta^*\cap N)\setminus\{0\} }a_n
\chi^{ n-\sum_{j=0}^k \min\langle \Delta_j,n\rangle r_j^* } $$ and
$\chi^{r_i^*}-\chi^{r_0^*}$, for $ i=1,\dots,k$. By
\cite[pp.~162-163]{a} or Lemma~2.2 in \cite{m3}, we already know
that
 the kernel of the surjective ring homomorphism
$$ \CC[\bar\sigma^\ve\cap \bar N]\longrightarrow\CC[\sigma^\ve\cap
(N\oplus\ZZ)] $$  is an ideal generated by the regular sequence
$\chi^{r_i^*}-\chi^{r_0^*}$, for $ i=1,\dots,k$. Therefore, it
suffices to show that any preimage of $f$ by this homorphism is in
the ideal generated by  $\bar f$ and $\chi^{r_i^*}-\chi^{r_0^*}$,
for $ i=1,\dots,k$.

For $n\in (\Delta^*\cap N)\setminus\{0\}$, the preimage of
$\chi^{(n,1)}$ by the ring homomorphism induced by    $
 N\oplus\ZZ^{k+1}\rightarrow N\oplus\ZZ$,
$(n,\alpha_0,\dots,\alpha_k)\mapsto(n,\alpha_0+\cdots+\alpha_k)$,
is a linear combination of  $\chi^{ n +\sum_{i=0}^k\alpha_jr_j^*}$
with $\sum_{j=0}^k\alpha_j=1$ and
$n+\sum_{j=0}^k\alpha_jr_j^*\in\bar\sigma^\ve$. But the last
condition means $\min\langle
n+\sum_{j=0}^k\alpha_jr_j^*,\Delta_l+r_l\rangle\ge0$, whence
$\alpha_l\ge -\min\langle n, \Delta_l \rangle $ for all $l$. Since
$$1=\sum_{j=0}^k\alpha_j\ge -\sum_{j=0}^k\min\langle n, \Delta_j \rangle=
-\min\langle n,\sum_{j=0}^k \Delta_j \rangle=-\min\langle n,
\Delta\rangle=1,$$ we get $\alpha_j=-\min\langle n, \Delta_j
\rangle$. It is also clear that a preimage of $\chi^{(0,1)}$
coincides with $\chi^{r_0^*}$ modulo  $\chi^{r_i^*}-\chi^{r_0^*}$,
for $ i=1,\dots,k$.
\end{pf}

Translating the above statement by the correspondence
(\ref{e:isom}) into homogeneous coordinates we get:

\begin{thm}\label{t:degener} Let $Y_{\Delta^*}\subset X_{\Delta^*}$ be a Calabi-Yau hypersurface given by
the equation $$ \sum_{n\in\Delta^*\cap N}a_n x^n
\prod_{\rho\in\Sigma_{\Delta^*}(1)}x_\rho=0,$$ where $a_n\in\CC$.
Then the image of $ Y_{\Delta^*} $ under the embedding $
X_{\Delta^*}\hookrightarrow
 X_{ {(\tilde\Delta_0\uplus\cdots\uplus\tilde\Delta_k)^*}}$   is a nef Calabi-Yau complete
intersection given by the equations $$
 \sum_{n \in\Delta^*\cap N }a_n x^{
n-\sum_{j=1}^k \min\langle \Delta_j,n\rangle e_j^* } \prod_{
 v_\rho\in\Delta_0+e_0 } x_\rho=0, \quad \prod_{
 v_\rho\in\Delta_i+e_i } x_\rho -\prod_{   v_\rho\in \Delta_0+e_0} x_\rho
=0,$$  for $i=1,\dots,k.$
\end{thm}

The ample Calabi-Yau hypersuface $Y_{\Delta^*}\subset
X_{\Delta^*}$ deforms to a generic nef Calabi-Yau complete
intersection $Y_{\hat\nabla_0,\dots,\hat\nabla_k}$ in the Fano
toric variety $ X_{
{(\tilde\Delta_0\uplus\cdots\uplus\tilde\Delta_k)^*}}=
X_{\hat\nabla_0+\cdots+\hat\nabla_k}$ corresponding to the
nef-partition $\hat\nabla_0+ \cdots+\hat\nabla_k$:
$$\Biggl(\sum_{j=1}^k a_{i,j} x^{ e_j^*-\delta_ie_i^*}+
 \sum_{n\in\Delta^* \cap N }a_{i,n}
x^{  n-\delta_ie_i^*-\sum_{j=1}^k \min\langle \Delta_j,n\rangle
e_j^* }\Biggr)\prod_{
 v_\rho\in\Delta_i+e_i } x_\rho=0 $$ for
$i=0,\dots,k$, where  $a_{i,j}, a_{i,n}\in\CC$ are the
coefficients, and $\delta_i=1$, if $i\ne0$, $\delta_0=0$. (Note
that the lattice points corresponding to the monomials are
precisely the lattice points of the polytope $\hat\nabla_i$ in
Proposition~\ref{p:lat}.)

\section{Degenerations and mirror contractions of Calabi-Yau complete
intersections.}\label{s:mirror}

In the previous section, we obtained a deformation of an ample
Calabi-Yau hypersurface $Y_{\Delta^*}\subset X_{\Delta^*}$ in a
Fano toric variety to a generic Calabi-Yau complete intersection
in the Fano toric variety $X_{\hat\nabla_0+\cdots+\hat\nabla_k}$.
Equivalently, we have a degeneration of a generic Calabi-Yau
complete intersection in $X_{\hat\nabla_0+\cdots+\hat\nabla_k}$ to
a generic Calabi-Yau hypersurface $Y_{\Delta^*}\subset
X_{\Delta^*}$. Let $\Sigma'_{\hat\nabla_0+\cdots+\hat\nabla_k}$ be
a maximal projective subdivision of the normal fan of the
reflexive polytope $\hat\nabla_0+\cdots+\hat\nabla_k$, and let
$\Sigma'_{\Delta^*}$ be the fan obtained by intersecting the cones
of $\Sigma'_{\hat\nabla_0+\cdots+\hat\nabla_k}$ with the linear
subspace  $M_\RR\subset\tilde  M_\RR\oplus
 \RR^k$.
Then we have a commutative diagram:
\begin{equation*}
\begin{array}{cccccc}
Y'_{\Delta^*}&\subset &X_{\Sigma'_{\Delta^*}}&\hookrightarrow &
X_{\Sigma'_{\hat\nabla_0+\cdots+\hat\nabla_k}}
\\
\downarrow& &\downarrow&&\downarrow& \\
Y_{\Delta^*}&\subset& X_{\Delta^*}&\hookrightarrow &
 X_{\hat\nabla_0+\cdots+\hat\nabla_k},
\end{array}
\end{equation*}
where $Y'_{\Delta^*}$ is a  crepant partial resolution of the
ample Calabi-Yau hypersurface $Y_{\Delta^*}$. The hypersurface
$Y'_{\Delta^*}$ deforms to a  Calabi-Yau complete intersection in
$X_{\Sigma'_{\hat\nabla_0+\cdots+\hat\nabla_k}}$. Correspondingly,
a generic Calabi-Yau complete intersection
$Y'_{\hat\nabla_0,\dots,\hat\nabla_k}$ in $
X_{\Sigma'_{\hat\nabla_0+\cdots+\hat\nabla_k}}$ degenerates to a
generic Calabi-Yau hypersurface $Y'_{\Delta^*}$ in
$X_{\Sigma'_{\Delta^*}}$. Now, if $\Sigma''_{\Delta^*}$ is a
maximal projective subdivision of the normal fan of the reflexive
polytope $\Delta^*$, which refines the fan $\Sigma'_{\Delta^*}$,
then we obtain a geometric transition  from a minimal Calabi-Yau
hypersurface to a minimal Calabi-Yau complete intersection:
\begin{equation}\label{e:trans}Y'_{\hat\nabla_0,\dots,\hat\nabla_k}\rightsquigarrow
Y'_{\Delta^*}\leftarrow Y''_{\Delta^*}, \end{equation}
 where $Y''_{\Delta^*}$ is a maximal  projective crepant  partial
 resolution of the ample
Calabi-Yau hypersurface $Y_{\Delta^*}\subset X_{\Delta^*}$.

By  Morrison's conjecture  in \cite{mor}, every geometric
transition between Calabi-Yau manifolds should correspond to a
mirror geometric transition between the mirror partners of the
original Calabi-Yau manifolds with the roles of degeneration and
contraction reversed. By the Batyrev-Borisov mirror symmetry
construction in \cite{b} and \cite{bb} we know that the mirror of
the Calabi-Yau hypersurface $Y''_{\Delta^*}$ is a nondegenerate
Calabi-Yau hypersurface in a maximal  projective crepant  partial
resolution of the Fano toric variety $X_\Delta$ and the mirror of
the Calabi-Yau complete intersection
$Y'_{\hat\nabla_0,\dots,\hat\nabla_k}$ is a nondegenerate
Calabi-Yau complete intersection in a maximal  projective crepant
partial resolution of the Fano toric variety
$X_{\tilde\Delta_0+\cdots+\tilde\Delta_k}$, corresponding to the
nef-partition $\tilde\Delta_0+\cdots+\tilde\Delta_k$ dual to
$\hat\nabla_0+\cdots+\hat\nabla_k$. We will explicitly construct a
natural  geometric transition between the mirror Calabi-Yau
varieties, which we expect to be the mirror of (\ref{e:trans}).

By Proposition~\ref{p:1iso}, the dual of the reflexive polytope
$\tilde\Delta_0+\cdots+\tilde\Delta_k$ is the convex hull
 $\hat\nabla_0\uplus\dots\uplus\hat\nabla_k$.
 By the construction,  the image of
 $\hat\nabla_0\uplus\dots\uplus\hat\nabla_k$ under the natural
 projection $\tilde N_\RR=N_\RR\oplus\RR^k\rightarrow N_\RR$ is
 the reflexive polytope $\Delta^*$.
It is not difficult to see that a subdivision of  the normal fan
of $\Delta$ lifts to a subdivision $\Sigma'_{
\tilde\Delta_0+\cdots+\tilde\Delta_k}$ of the normal fan of
$\tilde\Delta_0+\cdots+\tilde\Delta_k$. Without loss of
generality, we can assume that $\Sigma'_{
\tilde\Delta_0+\cdots+\tilde\Delta_k}$  is a maximal projective
subdivision of the normal fan of
$\tilde\Delta_0+\cdots+\tilde\Delta_k$ which maps to a
 maximal projective subdivision $\Sigma'_\Delta$  of the normal fan of $\Delta$ under the projection $\tilde N_\RR\rightarrow N_\RR$.
Hence, we get a toric morphism $X_{\Sigma'_{
\tilde\Delta_0+\cdots+\tilde\Delta_k}}\rightarrow
X_{\Sigma'_\Delta}$. The next result explicitly finds the image
under this morphism of a generic Calabi-Yau complete intersection
$Y'_{\tilde\Delta_0,\dots,\tilde\Delta_k}$ in $X_{\Sigma'_{
\tilde\Delta_0+\cdots+\tilde\Delta_k}}$ corresponding to the
nef-partition $\tilde\Delta_0+\cdots+\tilde\Delta_k$.

\begin{thm} Let $Y'_{\tilde\Delta_0,\dots,\tilde\Delta_k} \subset X_{\Sigma'_{
\tilde\Delta_0+\cdots+\tilde\Delta_k}}$ be a Calabi-Yau complete
intersection   given by the equations $$
\biggl(1-\sum_{m\in\Delta_i\cap M} a_m x^{m+e_i}\biggr)\prod_{
v_\rho\in  \hat\nabla_i} x_\rho=0,\quad i=0,\dots,k.$$ where
$a_m\in\CC$.  Then the image of
$Y'_{\tilde\Delta_0,\dots,\tilde\Delta_k}$ under the contraction
$X_{\Sigma'_{ \tilde\Delta_0+\cdots+\tilde\Delta_k}}\rightarrow
X_{\Sigma'_\Delta}$   is a nef Calabi-Yau   hypersurface given by
the equation \begin{equation}\label{e:image}
 \Biggl(1-\prod_{i=0}^k\Biggl(\sum_{m\in\Delta_i\cap M} a_m x^{m}\Biggr)\Biggr)\prod_{
v_\rho\in \Delta^* } x_\rho=0.\end{equation}
\end{thm}

\begin{pf} Note that the intersection of the Calabi-Yau complete
intersection $Y'_{\tilde\Delta_0,\dots,\tilde\Delta_k}$ with the
dense affine torus $\TT=\spe(\CC[M\oplus\ZZ^k])$ is a complete
intersection  given by the equations $1-\sum_{m\in\Delta_i\cap M}
a_m \chi^{m+e_i}=0$, for $i=0,\dots,k$. We can find the image of
$Y'_{\tilde\Delta_0,\dots,\tilde\Delta_k}$ as the closure of the
image of the affine complete intersection by the projection of
tori $\spe(\CC[M\oplus\ZZ^k])\rightarrow \spe(\CC[M]) $ induced by
the injective lattice homomorphism $M\subset M\oplus\ZZ^k$. By
eliminating  $\chi^{e_i}$, for $i=1,\dots,k,$ and
$\chi^{e_0}=\prod_{i=1}^k\chi^{-e_i}$ from the above equations we
get the equation
$$1-\prod_{i=0}^k \Biggl(\sum_{m\in\Delta_i\cap
M} a_m \chi^{m} \Biggr)=0$$
of the image in the affine torus
$\spe(\CC[M]) $. The Zariski closure of this affine hypersurface
is   the  nef Calabi-Yau hypersurface given by the equation
(\ref{e:image}) in homogeneous coordinates of
$X_{\Sigma'_\Delta}$.
\end{pf}

Denote by  $Y'_\Delta $ the    nef Calabi-Yau   hypersurface in
$X_{\Sigma'_\Delta}$ given by the equation (\ref{e:image}). Notice
that such a hypersurface is not generic. The geometric transition
from a generic Calabi-Yau complete intersection
$Y'_{\tilde\Delta_0,\dots,\tilde\Delta_k} \subset X_{\Sigma'_{
\tilde\Delta_0+\cdots+\tilde\Delta_k}}$ is completed by a
smoothing of $Y'_\Delta $ to a nondegenerate Calabi-Yau
hypersurface $Y''_\Delta$ in $X_{\Sigma'_\Delta}$:
\begin{equation}\label{e:mirtrans} Y''_\Delta\rightsquigarrow
Y'_{\Delta}\leftarrow  Y'_{\tilde\Delta_0,\dots,\tilde\Delta_k},
\end{equation}

As in \cite{mor}, we expect that the   natural geometric
transitions (\ref{e:trans}) and (\ref{e:mirtrans}) are mirror to
each other, which can be explained by a natural mirror map between
the complex and K\"ahler moduli spaces. In particular, the
degeneration $Y'_{\hat\nabla_0,\dots,\hat\nabla_k}\rightsquigarrow
Y'_{\Delta^*}$ and the mirror contraction $
Y'_{\tilde\Delta_0,\dots,\tilde\Delta_k}\rightarrow Y'_{\Delta}$
should correspond to certain parts of the respective compactified
complex and K\"ahler moduli spaces of the mirror pair of
Calabi-Yau complete intersections
$Y'_{\hat\nabla_0,\dots,\hat\nabla_k}$ and
$Y'_{\tilde\Delta_0,\dots,\tilde\Delta_k}$.

\section{Degenerations of the main periods and hypergeometric
series.}\label{s:degen}

We will support the mirror correspondence of the geometric
transitions by showing   that the degeneration of the main periods
(determining the mirror map) for Calabi-Yau complete intersection
$Y'_{\hat\nabla_0,\dots,\hat\nabla_k}$ and the hypersurface
$Y''_\Delta$ coincide with the main periods of the minimal
Calabi-Yau $Y''_{\Delta^*}$ and
$Y'_{\tilde\Delta_0,\dots,\tilde\Delta_k}$, respectively.

First, we recall the definition of the main period  for the
nondegenerate Calabi-Yau  hypersurface $Y''_\Delta$ from \cite{b1}
(also, see \cite[Sec.~6.3.4]{ck}. Fix an integer basis
$u_1,\dots,u_d$ for the lattice $M$. Then
$t_j=\prod_{\rho\in\Sigma_\Delta'(1)} x_\rho^{\langle
u_j,v_\rho\rangle}$, for $j=1,\dots,d$, are the coordinates on the
dense torus $\TT_N=\spe(\CC[M])=N\otimes_\ZZ\CC^*\subset
X_{\Sigma_\Delta'}$. Let $f_\Delta=1-  \sum_{m\in\partial\Delta
\cap M} b_m t^{m}$ be the Laurent polynomial  determining the
hypersurface $Y''_\Delta\cap\TT_N$, and let $\gamma\subset\TT_N$
be the cycle defined by $|t_1|=\cdots=|t_d|=1$, then the main
period for
 $Y''_\Delta$ equals the Euler integral
 $$\Phi_{Y''_\Delta}(\beta)=
 \frac{1}{(2\pi\sqrt{-1} )^d}\int_\gamma\frac{1}{f_{\Delta}}\frac{dt_1}{t_1}\wedge\cdots\wedge
\frac{dt_d}{t_d},$$ where $\beta=(b_m)\in\CC^{l(\Delta)-1}$. The
function $\Phi_{Y''_\Delta}(\beta)$ is called hypergeometric since
it  satisfies the GKZ hypergeometric system of differential
equations (see \cite[Theorem 14.2]{b1}). It can be found as a
power series expansion in the variables $b_m$:
$$\Phi_{Y''_\Delta}(\beta)=\sum_{l\in L_\Delta} (\sum_{m\in\partial\Delta
\cap M}l_m)!\prod_{m\in\partial\Delta \cap
M}\frac{b_m^{l_m}}{l_m!},$$ where
$L_\Delta=\{(l_m)_{m\in\partial\Delta \cap M}\mid \sum_{m
\in\partial\Delta \cap M} l_m m=0,
l_m\in\ZZ_{\ge0}\,\forall\,m\}$. It can also be written in terms
of the local coordinates on the complex moduli of $Y''_{\Delta}$
at a maximally unipotent boundary point (see \cite{bs,ck}). As
$Y''_{\Delta}$ degenerates to $Y'_\Delta$ the  hypergeometric
function  $\Phi_{Y''_\Delta}(\beta)$ will degenerate to the Euler
integral
$$\Phi_{Y'_\Delta}(\alpha)=
 \frac{1}{(2\pi\sqrt{-1} )^d}\int_\gamma\frac{1}{g_{\Delta}}\frac{dt_1}{t_1}\wedge\cdots\wedge
\frac{dt_d}{t_d},$$ where    $g_{\Delta}= 1-\prod_{i=0}^k
(\sum_{m\in\Delta_i\cap M} a_m t^{m} )$ determines the
hypersurface $Y'_\Delta\cap\TT_N$, and
$\alpha=(a_m)\in\CC^{l(\Delta_0)}\times\cdots\times\CC^{l(\Delta_k)}$.
Similar to \cite[Ex.~14.5]{b1}, substituting
$$\frac{1}{g_{\Delta}}=\sum_{j=0}^\infty \prod_{i=0}^k\Biggl(\sum_{m\in\Delta_i\cap M} a_m t^{m} \Biggr)^j$$
into the above integral and applying the Cauchy residue theorem
gives
\begin{equation}\label{e:hyperg}\Phi_{Y'_\Delta}(\alpha)=\sum_{l\in
L_{\tilde\Delta_0,\dots,\tilde\Delta_k}}\prod_{i=0}^k (\sum_{m\in
\Delta_i \cap M}l_m)!\prod_{m\in \Delta_i \cap
M}\frac{a_m^{l_m}}{l_m!},\end{equation} where
$$L_{\tilde\Delta_0,\dots,\tilde\Delta_k}=\Biggl\{(l_m) \mid \sum_{i=0}^k\sum_{m \in
\Delta_i \cap M} l_m m=0, \sum_{m \in \Delta_i \cap M} l_m=\sum_{m
\in \Delta_0 \cap M} l_m \,\forall\, i \Biggr\}$$ is a
subsemigroup of
$\ZZ_{\ge0}^{l(\Delta_0)}\oplus\cdots\oplus\ZZ_{\ge0}^{l(\Delta_k
)}$. But (\ref{e:hyperg}) is precisely the hypergeometric series
in \cite[Def.~6.1.1, Pr.~6.1.4]{bs} for the Calabi-Yau complete
intersection $ Y'_{\tilde\Delta_0,\dots,\tilde\Delta_k} $ given by
the Euler integral
$$\Phi_{Y'_{\tilde\Delta_0,\dots,\tilde\Delta_k}}(\alpha)= \frac{1}{(2\pi\sqrt{-1} )^{d+k}}
\int_{\tilde\gamma}\frac{1}{f_{\tilde\Delta_0}\cdots
f_{\tilde\Delta_k}}\frac{dt_1}{t_1}\wedge\cdots \wedge
\frac{dt_{d+k}}{t_{d+k}} ,$$ where
$t_1,\dots,t_d,t_{d+1},\dots,t_{d+k}$ are the coordinates on the
torus $\TT_{\tilde N}=\spe(\CC[\tilde M])$ corresponding to the
lattice basis $\{u_1,\dots,u_d,e_1,\dots,e_k\}$ of $\tilde
M=M\oplus\ZZ^k$, the cycle $\tilde\gamma$ is given by $|t_j|=1$
for
 $j=1,\dots,d+k$, and
$$f_{\tilde\Delta_0}= 1-\sum_{m\in\Delta_0\cap M} a_m t^{m
}\prod_{i=1}^k t^{-1}_{d+i},\quad f_{\tilde\Delta_i}=
1-\sum_{m\in\Delta_i\cap M} a_m t^{m }t_{d+i},\quad i=1,\dots,k,$$
determine the affine complete intersection $
Y'_{\tilde\Delta_0,\dots,\tilde\Delta_k}\cap \TT_{\tilde N}$. The
  series $\Phi_{Y''_\Delta}(\beta)$ and
$\Phi_{Y'_\Delta}(\alpha)$  are invariant under the natural torus
action $\TT_N$ on  the space of Laurent polynomials in the
variables $t_j$ and can be expressed in the local coordinates on
the complex moduli of  $Y''_\Delta$ as in \cite[Sec.~6.3.4]{ck}
Therefore,
  the main period
$\Phi_{Y''_\Delta}(\beta)$ of a minimal Calabi-Yau hypersurface
$Y''_\Delta$  degenerates to the main period
$\Phi_{Y'_{\tilde\Delta_0,\dots,\tilde\Delta_k}}(\alpha)$ of a
minimal Calabi-Yau complete intersection $
Y'_{\tilde\Delta_0,\dots,\tilde\Delta_k} $ as
$Y''_\Delta\rightsquigarrow Y'_{\Delta}$. The main period
$\Phi_{Y'_{\tilde\Delta_0,\dots,\tilde\Delta_k}}(\alpha)$
determines the mirror map between the complex moduli of
$Y'_{\tilde\Delta_0,\dots,\tilde\Delta_k}$ and the K\"ahler moduli
of the mirror partner $Y'_{\hat\nabla_0,\dots,\hat\nabla_k}$,
which allows to compute the instanton numbers of rational curves
in the latest Calabi-Yau as explained in \cite{ck}.

Similarly to the above, let us consider the main period for the
Calabi-Yau complete intersection
$Y'_{\hat\nabla_0,\dots,\hat\nabla_k}$:
$$\Phi_{_{\hat\nabla_0,\dots,\hat\nabla_k}}(\alpha)= \frac{1}{(2\pi\sqrt{-1} )^{d+k}}
\int_{\tilde\gamma}
\prod_{i=0}^k\frac{a_{i,i}}{f_{\hat\nabla_i}}\frac{dt_1}{t_1}\wedge\cdots
\wedge \frac{dt_{d+k}}{t_{d+k}} ,$$ where the coordinates on the
torus $\TT_{\tilde M}=\spe(\CC[\tilde N])$ corresponding to the
dual basis $\{u^*_1,\dots,u^*_d,e^*_1,\dots,e^*_k\}$ of $\tilde N$
are denoted by $t_1,\dots,t_{d+k}$ again (abusing  notation), and
$$f_{\hat\nabla_i}=t_{d+i}^{-1}\Biggl(a_{i,0}+\sum_{j=1}^k a_{i,j} t_{d+j}+
 \sum_{n\in\partial\Delta^* \cap N }a_{i,n}
t^n\prod_{j=1}^kt_{d+j}^{-\min\langle \Delta_j,n\rangle }\Biggr),
\quad i=0,\dots,k,$$ (with $t_{d+i}=1$ if $i=0$) determine the
affine complete intersection
$Y'_{\hat\nabla_0,\dots,\hat\nabla_k}\cap\TT_{\tilde M}$, and
$\alpha=(a_{*,*})\in\CC^{(k+1)(l(\Delta^*)+k)}$. Rewriting $
a_{i,i}/f_{\hat\nabla_i}$ as a series $$\sum_{s=0}^\infty
  \Biggl( \sum_{j\ne i}\frac{ a_{i,j}}{-a_{i,i}}t_{d+i}^{-1}t_{d+j}+
 \sum_{n\in\partial\Delta^* \cap N }\frac{a_{i,n}}{-a_{i,i}}t_{d+i}^{-1}
t^n\prod_{j=1}^kt_{d+j}^{-\min\langle \Delta_j,n\rangle }
\Biggr)^s$$ (here $j\ne i$ means
$j\in\{0,\dots,k\}\setminus\{i\}$) and applying the Cauchy residue
formula similar to \cite[Ex.~14.5]{b1} gives
$$\Phi_{_{\hat\nabla_0,\dots,\hat\nabla_k}}(\alpha) =\sum_{l\in
L_{\hat\nabla_0,\dots,\hat\nabla_k}}\prod_{i=0}^k
\frac{l_{i}!}{(-a_{i,i})^{l_{i}}}\prod_{j\ne i}
\frac{a_{i,j}^{l_{i,j}}}{l_{i,j}!}\prod_{n\in\partial\Delta^* \cap
N}\frac{a_{i,n}^{l_{i,n}}}{l_{i,n}!},$$ where $l_i=\sum_{j\ne i}
l_{i,j}+\sum_{n\in\partial\Delta^* }l_{i,n}$, and
 $L_{\hat\nabla_0,\dots,\hat\nabla_k}$ consists   of vectors $(l_{*,*})\in\ZZ_{\ge0}^{(k+1)(l(\Delta^*)+ k-1)} $
 such that
 $$\sum_{i=0}^k\Biggl(\sum_{j\ne i}l_{i,j}(\delta_je_j^*-\delta_ie_i^*)+\sum_{n\in\partial\Delta^*  }
 l_{i,n} \biggr(n-\delta_ie_i^*-\sum_{j=1}^k \min\langle \Delta_j,n\rangle e_j^* \biggl)\Biggr)=0,$$
where $\delta_j=1$ if $j\ne0$ and $\delta_0=0$.
 From Theorem~\ref{t:degener}, degeneration
$Y'_{\hat\nabla_0,\dots,\hat\nabla_k}\rightsquigarrow
Y'_{\Delta^*}$ corresponds to setting  $a_{0,0}=a_0$, $a_{0,j}=0$,
$a_{0,n}=a_n$,
   $a_{i,i}=1$, $a_{i,0}=-1$
   $a_{i,i}=1$, $a_{i,0}=-1$, $a_{i,j}=0$, $a_{i,n}=0$ for $1\le i,j\le k$, $i\ne j$,  and $n\in\partial\Delta^* \cap N$.
Coefficients $a_{*,*}=0$ force the vanishing of the  terms in the
  series $\Phi_{_{\hat\nabla_0,\dots,\hat\nabla_k}}(\alpha)$
unless the corresponding $l_{*,*}$ equals zero. Hence, the
nonvanishing  terms in the series correspond to
$l_0=\sum_{n\in\partial\Delta^*\cap N }l_{0,n}$, $l_i= l_{i,0}$,
for $i=1,\dots,k$, and
$$\sum_{n\in\partial\Delta^* \cap N }
 l_{0,n} \biggr(n -\sum_{j=1}^k \min\langle \Delta_j,n\rangle e_j^* \biggl)+
 \sum_{i=1}^k l_{i,0} (- e_i^*)=0,$$
and the series $\Phi_{_{\hat\nabla_0,\dots,\hat\nabla_k}}(\alpha)$
degenerates to $$\sum_{l}\Biggl(\prod_{i=1}^k
\frac{l_{i}!}{(-1)^{l_{i}}}
\frac{(-1)^{l_{i,0}}}{l_{i,0}!}\Biggr)\frac{l_{0}!}{(-a_0)^{l_{0}}}\prod_{n\in\partial\Delta^*
\cap N}\frac{a_{n}^{l_{0,n}}}{l_{0,n}!}=\sum_{l\in L_{\Delta^*}}
\frac{l_{0}!}{(-a_0)^{l_{0}}}\prod_{n\in\partial\Delta^* \cap
N}\frac{a_{n}^{l_{n}}}{l_{n}!},$$ where
$$L_{\Delta^*}=\biggl\{(l_n)_{n\in\partial\Delta^* \cap N}\mid \sum_{n
\in\partial\Delta^* \cap N} l_n n=0\biggr\}\subset
\ZZ_{\ge0}^{l(\Delta^*)-1}, \quad
l_0=\sum_{n\in\partial\Delta^*\cap N }l_{ n}.$$  But the last
series is  the hypergeometric series corresponding to the main
period of Calabi-Yau hypersurface $Y''_{\Delta^*}$ at the
maximally unipotent boundary point (see \cite{b1,ck}). Thus, we
showed that  the main period
$\Phi_{_{\hat\nabla_0,\dots,\hat\nabla_k}}(\alpha)$ of a minimal
Calabi-Yau complete intersection $
Y'_{\hat\nabla_0,\dots,\hat\nabla_k}$ degenerates to the main
period $\Phi_{Y''_\Delta}(\alpha_0)$ of a minimal Calabi-Yau
hypersurface $Y''_{\Delta^*}$ as
$Y'_{\hat\nabla_0,\dots,\hat\nabla_k}\rightsquigarrow
Y'_{\Delta^*}$.

There is more work to be done in  computing the mirror map itself
from the main periods of Calabi-Yau, but we will finish this paper
by explaining the expected relationship of the complex and
K\"ahler moduli of the Calabi-Yau varieties involved in our
geometric transitions. For definitions and notation we refer to
the book \cite{ck}.

A minimal Calabi-Yau variety $V$ has a K\"ahler cone $K(V)$,  a
complexified  K\"ahler space
$$K_\CC(V)=\{\omega\in H^2(V,\CC)\mid {\rm Im}(\omega)\in
K(V)\}/\im H^2(V,\ZZ),$$ and the complexified K\"ahler moduli
space ${\cal K}(V)=K_\CC(V)/{\rm Aut}(V)$. For a nondegenerate
Calabi-Yau complete intersection $V$ in a maximal projective
partial crepant resolution $\xs$ of a Gorenstein Fano toric
variety, one considers a {\it toric  part of the K\"ahler cone}
$K(V)_{\rm toric}=K(V)\cap H^2_{\rm toric}(V)$, where $H^2_{\rm
toric}(V)$ is the image of the restriction map
$H^2(X_\Sigma)\rightarrow H^2(V)$, and the corresponding {\it
toric  K\"ahler moduli space} ${\cal K}(V)_{\rm toric}$.

On the complex side, one also considers a part of the complex
moduli space of  a Calabi-Yau complete intersection $V\subset\xs$.
Let $V$ be
  the closure   of the affine complete intersection
$$\sum_{m\in\Delta_i\cap \ZZ^d} a_{i,m} t^m=0, \quad i=0,\dots,k$$
in  $(\CC^*)^d\subset  \xs$, where
$\Delta=\Delta_0+\cdots+\Delta_k$ is a Minkowski sum decomposition
of a reflexive polytope. Then the {\it polynomial moduli space} of
the complete intersection $V$ can be constructed similar to
\cite[Sect.~13]{bc} as a geometric quotient ${\cal M}(V)_{\rm
poly}=U/{\rm Aut}(\xs)$, where $U$ is an open subset in
$\PP(L(\Delta_1\cap \ZZ^d))\times\cdots\times\PP(L(\Delta_k\cap
\ZZ^d))$ corresponding to a subset of the set of  quasismooth
complete intersections (see \cite{m0}) with $L(\Delta_i\cap
\ZZ^d)$ denoting the vector space of Laurent polynomials
$\sum_{m\in\Delta_i\cap \ZZ^d} a_{i,m} t^m$.

In practice, one replaces the toric   K\"ahler moduli space and
the polynomial moduli space of a Calabi-Yau intersection $V$ with
suitable compactifications  realized as toric varieties associated
with a secondary fan (see \cite{ck}). We will denote them by
$\overline{\cal K}(V)_{\rm toric}$ and $\overline{\cal M}(V)_{\rm
poly}$, respectively.

The geometric transitions
$Y'_{\hat\nabla_0,\dots,\hat\nabla_k}\rightsquigarrow
Y'_{\Delta^*}\leftarrow Y''_{\Delta^*}$
 and $Y''_\Delta\rightsquigarrow
Y'_{\Delta}\leftarrow  Y'_{\tilde\Delta_0,\dots,\tilde\Delta_k}$
induce inclusions $$\overline{\cal M}( Y''_{\Delta^*})_{\rm
poly}\subset \overline{\cal M}(
Y'_{\hat\nabla_0,\dots,\hat\nabla_k})_{\rm poly},\quad
\overline{\cal M}( Y'_{\tilde\Delta_0,\dots,\tilde\Delta_k})_{\rm
poly}\subset \overline{\cal M}( Y''_\Delta)_{\rm poly},$$
corresponding to degenerations on the complex side, and inclusions
$$\overline{\cal K}(Y'_{\hat\nabla_0,\dots,\hat\nabla_k})_{\rm toric}\subset
\overline{\cal K}(Y''_{\Delta^*})_{\rm toric},\quad \overline{\cal
K}(Y''_\Delta)_{\rm toric}\subset \overline{\cal
K}(Y'_{\tilde\Delta_0,\dots,\tilde\Delta_k})_{\rm toric}, $$
corresponding to contractions on the K\"ahler side. These should
fit into the following commutative diagrams:
$$
\begin{array}{cccccccccccc}
&\overline{\cal M}( Y''_{\Delta^*})_{\rm poly}&\subset&
\overline{\cal M}( Y'_{\hat\nabla_0,\dots,\hat\nabla_k})_{\rm
poly}& &\overline{\cal M}(
Y'_{\tilde\Delta_0,\dots,\tilde\Delta_k})_{\rm
poly}&\subset& \overline{\cal M}( Y''_\Delta)_{\rm poly}&  \\
&\downarrow & &\downarrow & &\downarrow & &\downarrow & \\
&\overline{\cal K}(Y''_\Delta)_{\rm toric}&\subset& \overline{\cal
K}(Y'_{\tilde\Delta_0,\dots,\tilde\Delta_k})_{\rm toric},&
&\overline{\cal K}(Y'_{\hat\nabla_0,\dots,\hat\nabla_k})_{\rm
toric}&\subset&\overline{\cal K}(Y''_{\Delta^*})_{\rm toric},&
\end{array}
$$ where the vertical arrows are the mirror morphisms between the
complex and K\"ahler moduli spaces. Moreover,   the degenerations
of the main periods of Calabi-Yau varieties (calculated above)
should induce degenerations of the mirror morphisms between the
ambient moduli to the mirror morphisms of the enclosed moduli. The
inclusions of the  moduli can be described explicitly in terms of
  the inclusions of the respective secondary fans.

\end{document}